[1994/12/01]
\documentclass[draft]{article}
\title{The Radon-Nikodym problem for approximately proper equivalence relations}
\author{Jean Renault}
\date{November 23, 2002}

\usepackage{amsmath,amsthm,amscd}

\chardef\bslash=`\\ 

\newtheorem{thm}{Theorem}[section]
\newtheorem{cor}[thm]{Corollary}
\newtheorem{lem}[thm]{Lemma}
\newtheorem{prop}[thm]{Proposition}

\theoremstyle{definition}
\newtheorem{defn}{Definition}[section]
\newtheorem{ex}{Example}[section]

\theoremstyle{remark}
\newtheorem{rem}{Remark}[section]

\newcommand{\thmref}[1]{Theorem~\ref{#1}}

\newcommand{\lemref}[1]{Lemma~\ref{#1}}
\newcommand{\propref}[1]{Proposition~\ref{#1}}
\newcommand{\corref}[1]{Corollary~\ref{#1}}
\newcommand{\defnref}[1]{Definition~\ref{#1}}
\newcommand{\exref}[1]{Example~\ref{#1}}

\newcommand{\eval}[2][\right]{\relax
  \ifx#1\right\relax \left.\fi#2#1\rvert}

\def\projlim{\mathop{\oalign{lim\cr
\hidewidth$\longleftarrow$\hidewidth\cr}}}

\def\indlim{\mathop{\oalign{lim\cr
\hidewidth$\longrightarrow$\hidewidth\cr}}}

\begin{document}
\maketitle
\markboth{Jean Renault}
{AP equivalence relations}
\begin{abstract}
We study the Radon-Nikodym problem for approximately proper equivalence
relations and more specifically the uniqueness of certain Gibbs states. 
One of our
tools is a variant of the dimension
group introduced in the study of AF algebras. As
applications, we retrieve sufficient conditions for the
uniqueness of traces on AF algebras and parts of
the Perron-Frobenius-Ruelle theorem.\footnote{


Cocycles. Radon-Nikodym derivative. Dimension groups. Perron-Frobenius theorem.
\ }

\end{abstract} 
\renewcommand{\sectionmark}[1]{}

\section{Introduction}

The motivation of this work is a problem in the theory
of C$^*$-algebras, namely the study of the KMS
states of some automorphism groups of the Cuntz
algebras and their generalizations studied in \cite{el:kms}.
While the crucial role of the Perron-Frobenius theorem
in this problem has been noticed from the
beginning (see for example \cite{efw:kms, eva:markov}),
the application of Ruelle's version of this theorem
is more recent (\cite{kp:pressure, exe:RPF, ren:AF}). The existing
proofs of the Perron-Frobenius-Ruelle theorem, in
particular \cite{fan:Ruelle}, use heavily a sequence of
expectations associated with the asymptotic algebra. The
purpose of this work is to use the formalism of groupoids
(cf. \cite{ren:approach}), approximately proper
equivalence relations and dimension groups (cf. \cite{goo:poag}) to provide a
convenient setting for these proofs.

Given a groupoid
$G$ on a space $X$ and a cocycle $D\in
Z^1(G,{\bf R}_+^*)$, the Radon-Nikodym
problem mentioned in the title is the
study of the probability measures
on $X$ which are quasi-invariant with respect to $G$
(the definition is recalled in Section 2) and which
admit this cocycle as Radon-Nikodym derivative. 
When
$G=R$ is an approximately
proper (abbreviated as AP) equivalence relation, i.e. an increasing
union of proper equivalence relations $R_n$, the
cocycle defines a sequence of expectations
$E_n$ with range $C(X/R_n)$. The solutions of the
Radon-Nikodym problem are exactly the measures which
factor through $E_n$ for all $n$. A classical
example of this situation is provided by statistical
mechanics on a lattice
$\Lambda$. The sequence $(R_n)$ is defined by an
increasing sequence of finite sets $\Lambda_n$ with
union
$\Lambda$. The cocycle is the energy
cocycle, as in Section II.5 of \cite{ren:approach}. In this setting, the
solutions of the Radon-Nikodym problem are called Gibbs
states (this definition has been introduced by D.~Capocaccia in
\cite{cap:gibbs}). 

Our main concern is the uniqueness of
the solution of the Radon-Nikodym problem on an AP equivalence
relation. We give sufficient conditions for uniqueness
in two cases. First we consider
quasi-product cocycles on AF equivalence relations.
The data consist of a Bratteli diagram and a labeling
of its edges. Then a convenient condition on
this labeling (\corref{AF unique state},$(ii)$) guarantees uniqueness. Our
result covers (and was inspired by)
A.~T\"or\"ok's work \cite{toe:AF} (pointed to me by O.~Bratteli) on uniqueness
of traces on AF C$^*$-algebras.
The second case is the classical  setting of Ruelle's Perron Frobenius theorem and
gives the part of P.~Walters' Theorem 8 in \cite { wal:78}  concerning the
transpose of the Ruelle operator (\corref{Theorem 6} and \propref{Ruelle}). The
eigenvalue problem for this operator amounts to a Radon-Nikodym problem on the
semi-direct groupoid
$G(X,T)$. Our method, outlined in Section 4.1 of \cite{ren:AF}, is to solve 
first the Radon-Nikodym problem on the asymptotic equivalence relation $R(X,T)$
which is approximately proper. If it has a unique solution, this is also
a solution of the initial problem.

This work is organized as follows. The definition and some
examples of AP equivalence relations are given in the
first section. The second section deals with cocycles. We
are only concerned here with cocycles with values in the
multiplicative group ${\bf R}_+^*$ of
strictly positive real numbers. Let $D$ be a cocycle defined on the
AP equivalence relation $R=\cup R_n$ on the compact space
$X$. Then its restriction to $R_n$ can be written as a
coboundary. This provides a normalized potential $\rho_n$
and an expectation $E_n$ from $C(X)$ to $C(X/R_n)$, where
$C(X)$ denotes the space of real-valued continuous
functions on the compact space $X$. This is the sequence of expectations
mentioned earlier. A necessary and sufficient condition
for the unique ergodicity of $D$ is that for all $f\in
C(X)$, the variation of $E_n(f)$ tends to $0$. The theory
of dimension groups (this means here inductive limits of
$C(X_n)$'s, viewed as
ordered vector spaces, under positive linear maps), as
developped by K.~Goodearl in \cite{goo:poag}, is well suited to our
problem. Indeed, the solutions of the Radon-Nikodym
problem are the states of the dimension group associated
to the sequence $(E_{n,n-1})$ (where $E_n=E_{n,n-1}\circ
E_{n-1}$). Therefore, we recall in an appendix some
results of this theory. When the vector spaces $C(X_n)$
have finite dimension, there is an elementary sharp
estimate of the rate at which Markovian (i.e. unital)
operators contract the variation (\lemref{torok}) (this is essentially the
same estimate which is used by T\"or\"ok in \cite{toe:AF}; it is so natural that
it must have been noticed by other people). Rather surprisingly, this contraction
rate admits an easy estimate in terms of non-unital positive linear maps
(\lemref{key estimate}). This gives an elementary proof of the Perron-Frobenius
theorem for primitive matrices. This estimate is also used
in the third section devoted to quasi-product cocycles on
AF equivalence relations. This section is very close to
the sections 3 and 4 of \cite{eva:quasiproduct}, where
the emphasis is on measures rather than on cocycles. The
main result is \corref{AF unique state}, which gives a sufficient condition for
unique ergodicity as mentioned earlier. The section four introduces some
definitions such as $\pi$-cover, which are useful when 
dealing with local homeomorphism on  arbitrary compact
spaces. When $X$ is a compact space, we find more
convenient to use the entourages associated to a finite
open cover rather than those defined by a compatible
metric. This section also contains some estimates of the
variation. The section five studies the case of a
stationary system defined by a single surjective local
homeomorphism $T:X\rightarrow X$, where $X$ is compact.
The main result, \thmref{Walters}, is well-known: it gives the unique
ergodicity of the cocycle $D$ defined by a potential
$g\in C(X,{\bf R}_+^*)$ (or by a sequence of potentials
$(g_n)$) under the usual assumptions on the dynamical
system ($T$ is assumed to be expansive and
$R(X,T)$ to be minimal) and on the potentials (Walter's
condition). The proof is more elementary than most in the
sense that it does not use the Schauder-Tychonoff
theorem. However, the key step relies on the
Ascoli-Arzela theorem just as in \cite{wal:78}. It may be
that the estimate of the contraction rate of the
variation gives a more precise proof, with an estimate of
the speed of convergence; but this is not done here. In
the last section, the results are applied to the transfer
operator.
\vskip1cm

{\bf Acknowledgments} 
\vskip.3cm
A major part of this work was done
at the Centre for Advanced Study in Oslo while I held a
Visiting Fellowship under the program ''Non-commutative Phenomena in Mathematics
and Theoretical Physics''. I heartily thank its staff for an excellent stay and the
organizers, in particular M.~Landstad , for their invitation. I am also grateful
to E.~Alfsen, O.~Bratteli, R.~Exel, A.~Kumjian and C.~Skau for
fruitful discussions and suggestions.

\section{AP equivalence relations}

In this section, $X$ is a locally compact second countable Hausdorff space and $R$
is an equivalence relation on $X$. For the sake of simplicity, we only consider
equivalence relations with countable equivalence classes. We also
denote by
$R\subset X\times X$ its graph.

\begin{defn}\label{proper} The equivalence relation $R$ on $X$ will be called
{\it proper and \'etale} if its quotient space is Hausdorff and its quotient map
$X\rightarrow X/R$ is a local homeomorphism.
\end{defn}

Endowed with the product topology of $X\times X$, $R$ is a locally compact groupoid
which is \'etale. This simply means that the
projection maps
$r,s: R\rightarrow X$ are local homeomorphisms. Every open cover $\mathcal U$ of
$X$ by open sections of
$\pi$ gives a cover $\{(U\times V)\cap R\}, U,V\in{\mathcal U}\}$ of $R$ by open
bisections (i.e. simultaneous sections of the projection maps $r,s$). The
construction of
\cite{ren:approach} applies and gives the C$^*$-algebra $C^*(R)$. Moreover this
groupoid is proper in the sense of \cite{dr:amenable}.  It is
known, (see \cite{mw:continuous trace}) that $C^*(R)$ has continuous
trace and that it is Rieffel-Morita equivalent to $C_0(X/R)$. Since our study
is limited to \'etale equivalence relations and groupoids, we shall often omit the
word
\'etale and say ``proper'' rather than ``proper and \'etale''.

\begin{defn}\label{AP} The equivalence relation $R$ on $X$ will be called {\it
approximately proper}, abbreviated {\it AP}, if there exists a sequence
$$X_0\xrightarrow{\pi_{1,0}} X_1\rightarrow\ldots\rightarrow
X_{n-1}\xrightarrow{\pi_{n,n-1}}X_n
\rightarrow\ldots$$
where $X_0=X$ and for each $n\ge 1$, $X_n$ is a
Hausdorff space and $\pi_{n,n-1}$ is a surjective local
homeomorphism such that
$$R=\{(x,y)\in X\times X:\quad \exists n\in{\bf N}: \pi_n(x)=\pi_n(y)\}$$
where
$\pi_n=\pi_{n,n-1}\circ\pi_{n-1,n-2}\ldots\circ\pi_{1,0}$.
\end{defn}

In the context of the above definition, we let $R_n$ be the equivalence
relation on
$X$ defined by
$\pi_n$.

\begin{prop}\label{groupoid inductive limit} Let $R$ be the equivalence relation
defined by the sequence $(\pi_n)$ as above.
\begin{enumerate}
\item $(R_n)$ is an increasing sequence of subsets of $X\times X$ and $R=\cup
R_n$. 
\item For
$m\le n$, $R_m$ is a closed and open subset of $R_n$. 
\item Endowed with the inductive
limit topology, $R$ is an \'etale locally compact groupoid.
\end{enumerate}
\end{prop}

\begin{proof} The assertion $(i)$ is obvious. For the assertion $(ii)$, we
introduce the equivalence relation $S$ on $X_m$ defined by the map $\pi_{n,m}$.
The diagonal
$\Delta_{X_m}$ is closed and open in $S$. Therefore,
$R_m=(\pi_m\times\pi_m)^{-1}(\Delta_{X_m})$ is open and closed in
$R_n=(\pi_m\times\pi_m)^{-1}(S)$. For the assertion $(iii)$, according to
\cite{ren:approach}, it suffices to construct a cover of $R$ consisting of open
locally compact bisections. For each $n\in{\bf N}$, let ${\mathcal U}_n$ be a
cover of $X$. Then, the family $\{(U\times V)\cap R_n\}$, where $U,V\in{\mathcal
U}_n$, is an open cover of $R_n$ and the union of these covers is the desired cover
of
$R$.
\end{proof}

\begin{cor}\label{C$^*$ inductive limit} Let $R=\cup R_n$ be an AP equivalence
relation as above. Then, $C^*(R)$ is the inductive limit of the
$(C^*(R_n))$'s. More precisely, $C^*(R_n)$ can be identified with a sub
C$^*$-algebra; these subalgebras are increasing and their union is dense.
\end{cor}

On the other hand, for all $n\in{\bf N}$, $C(X_n)$ can be identified with the
subalgebra
$(\pi_n^*C(X_n))$ of $C(X)$ and this sequence of subalgebras is decreasing. Its
intersection is the fixed point subalgebra, as defined below.

\begin{defn} Let $(X,R)$ be an \'etale equivalence relation. The {\it fixed point
subalgebra} is the subalgebra of $C(X)$:
$$C(X)^R=\{f\in C(X):(x,y)\in R\Rightarrow f(x)=f(y)\}.$$
Its spectrum will be denoted by $X_\infty$. One says that $(X,R)$ is
{\it irreducible} if $C(X)^R$ consists only of constant functions.
\end{defn}

Thus, when $(X,R)$ is an AP equivalence relation defined by a sequence $(X_n)$,
the fixed point subalgebra 
$C(X)^R$  is the intersection of the
$\pi_n^*(C(X_n))$'s. The inclusion $C(X)^R\subset
\pi_n^*C(X_n)$ defines a continuous surjective map $\pi_{\infty,n}:X_n\rightarrow
X_\infty$. In particular, we write $\pi_\infty=\pi_{\infty,0}:X\rightarrow
X_\infty$. These maps identify $X_\infty$ to the inductive limit of the sequence
$(\pi_{n,n-1}:X_{n-1}\rightarrow X_n)$. In the irreducible case, this space is
reduced to a point  while the equivalence relation $R$ is non-trivial.

\begin{ex}{\em AF equivalence relations.} (See \cite{ren:AF} and
\cite{gps:affable}.) Let
$(V,E)$ be a Bratteli diagram. Recall that this means an oriented graph, where the
vertices are stacked on levels
$n=0,1,2,\ldots$ and the edges run from a vertex of level $n-1$ to a vertex of
level $n$. We denote by $V(n)$ the set of vertices of the level $n$ and by $E(n)$
the set of edges from level $n-1$ to level $n$. We assume that  every vertex
$v$ emits  finitely many, but at least one, edges and that every vertex on a
level $n\ge 1$ receives at least one edge. An infinite path is a sequence of
connected edges
$x=x_1x_2\ldots$, where $x_1$ starts from level $0$. The space $X$ of
infinite paths has a natural totally disconnected topology, with the cylinder sets
$Z(x_1x_2\ldots x_n)$ as a basis. We define similarly the space $X_n$ of infinite
paths starting from level $n$ and we have the obvious projection maps
$\pi_{n,n-1}:X_{n-1}\rightarrow X_n$. This defines an AP equivalence relation
$R$ on $X$ called the tail equivalence relation of the Bratteli diagram. In
the sequel, following \cite{gps:affable}, we shall use the notation
$(X=X(V,E),R=R(V,E))$ to designate this AP equivalence relation; we shall call it
the tail  equivalence relation of the Bratteli diagram $(V,E)$. An explicit
construction of matrix units in
$C^*(R)$ shows  that it is an AF C$^*$-algebra admitting $(V,E)$ as Bratteli
diagram. Conversely, it is shown in
\cite{gps:affable} (and in \cite{ren:AF} when the space $X$ is compact) that every
AP equivalence relation on a space
$X$ which is locally compact and totally disconnected is the tail equivalence
relation of a Bratteli diagram. Such an equivalence relation is called an AF
equivalence relation. Properties of AF equivalence relations are studied in
\cite{gps:affable}, where an equivalent definition is used.
\end{ex}
\begin{ex} {\em Stationary equivalence relations.} Let $X$ be a locally compact
space and
$T:X\rightarrow X$ a local homeomorphism which is onto. We form the stationary
sequence 
$$X\xrightarrow{T} X\xrightarrow{T} X\xrightarrow{T}\ldots$$
The corresponding AP equivalence relation
$$R=\{(x,y)\in X\times X:\quad\exists n\in {\bf N}: T^nx=T^ny\}$$
plays an essential role in the study of $T$ when it is large enough. This
example will be developped later.
\end{ex}

The definition of an AP equivalence relation makes an implicit reference to a
defining sequence
$$X_0\xrightarrow{\pi_{1,0}} X_1\rightarrow\ldots\rightarrow
X_{n-1}\xrightarrow{\pi_{n,n-1}}X_n
\rightarrow\ldots$$

\begin{defn}\label{equivalent defining sequences} We shall say that two sequences
$$\begin{array}{ccccccccc}
X_0&\rightarrow& X_1 &\rightarrow& \ldots & \rightarrow  &X_m
 &\rightarrow &\ldots \\
Y_0&\rightarrow& Y_1 &\rightarrow &\ldots & \rightarrow & Y_n
 &\rightarrow&\ldots
\end{array}$$
as in \defnref{AP} are equivalent if there are subsequences $(m_k)$ and $(n_k)$
and surjective local homeomorphisms $X_{m_k}\rightarrow Y_{n_k}$ and
$Y_{n_k}\rightarrow X_{m_{k+1}}$ making commutative diagrams.
\end{defn}

This is an
equivalence relation. An other way to obtain the same definition is to define
first the contraction of a sequence $(X_n)$: it is the new sequence $({\underline
X}_k)$ defined by a subsequence $(n_k)$; explicitly, ${\underline
X}_k=X_{n_k}$ and ${\underline\pi}_{k,k-1}=\pi_{n_k,n_{k-1}}$. Then, we say that
the original sequence is a dilation of the new sequence. Two sequences are
equivalent iff they admit contractions which have a common dilation.

\begin{prop}\label{equivalence of defining sequences} Let $X$ be a locally compact
space.
\begin{enumerate}
\item Equivalent defining sequences $(X_n)$ and $(Y_n)$ as above with $X_0=Y_0=X$
define the same equivalence relation $R$ on $X$ and the same topology on $R$.
\item Conversely, if $X$ is compact, two
sequences $(X_n)$ and $(Y_n)$ as above with $X_0=Y_0=X$ which define the
same equivalence relation $R$ are equivalent.
\end{enumerate} 
\end{prop}

\begin{proof}
Let us call $(R_n)$ [resp. $(S_n)$] the sequence of equivalence relations on $X$
defined by
$(X_n)$ [resp. $(Y_n)$]. If $(X_n)$ and $(Y_n)$ are equivalent and if we have
subsequences $(m_k)$ and $(n_k)$ as in the definition, we have the inclusions
$R_{m_k}\subset S_{n_k}\subset R_{m_{k+1}}$ for all $k$. We also know from
\propref{equivalent defining sequences} $(ii)$ that these inclusion maps are open.
Therefore the sequences
$(R_n)$ and $(S_n)$ have the same union $R$ and the inductive limit topology is
the same. Suppose now that $X$ is compact and that the sequences
$(R_n)$ and $(S_n)$ have the same union $R$. Let us endow $R$ with the inductive
limit topology of the $(S_n)$'s. Let us fix $m$. Since $R_m$ is closed in $X\times
X$, it is a compact subset of $R$. Since it is covered by the union of the open
sets $S_n$'s, it is contained in some $S_n$. Similarly, any $S_n$ is contained in
some $R_m$. Therefore one can construct subsequences $(m_k)$ and $(n_k)$ such that
$R_{m_k}\subset S_{n_k}\subset R_{m_{k+1}}$ for all $k$. These inclusions give the
desired maps $X_{m_k}\rightarrow Y_{n_k}$ and
$Y_{n_k}\rightarrow X_{m_{k+1}}$.

\end{proof}

\section{Cocycles}

We shall only consider cocycles with values in the group $A={\bf R}$ or
equivalently
$A={\bf R}_+^*$. If $R$ is an equivalence relation on $X$, a cocycle with values
in $A$ is a map
$c:R\rightarrow A$ satisfying $c(x,y)+c(y,z)=c(x,z)$ for all $(x,y),(y,z)\in R$. A
cocycle $c$ is a coboundary if there is a map
$b:X\rightarrow A$, called a {\it potential}, such that $c(x,y)=b(x)-b(y)$. Two
cocycles are cohomologous if their difference is a coboundary. In our setting, $R$
is a topological groupoid
$A$ is a topological group and we assume that $c$ and $b$ as above are
continuous. We denote by $Z^1(R,A)$ the set of continuous $A$-valued cocycles. We
shall use the convention $c\in Z^1(R,{\bf R})$ and $D\in Z^1(R,{\bf R}_+^*)$.

One defines similarly cocycles on arbitrary groupoids. Let $G$ be a locally compact
groupoid with a continuous Haar system (cf.
\cite{ren:approach}, Chapter 1). A measure $\mu$ on its unit space $X=G^{(0)}$ is
said to be quasi-invariant if the measures
$\mu\circ\lambda$ and $(\mu\circ\lambda)^{-1}$ on
$G$ are equivalent. Then, the Radon-Nikodym derivative $\displaystyle
D_\mu={d(\mu\circ\lambda)\over d(\mu\circ\lambda)^{-1}}$ satisfies the cocycle
identity almost everywhere. When $G$ is an \'etale groupoid, for example a proper
and
\'etale equivalence relation or an AP equivalence relation as above, it carries the
canonical continuous Haar system consisting of counting measures $\lambda^x$ on
the fibers $G^x=r^{-1}(x)$ of the range map $r$. A measure $\mu$ on $X=G^{(0)}$ is
quasi-invariant if and only if for all open bisections $S\subset G$, the measures
$\sigma(S)_*(\mu_{|s(S)})$ and $\mu_{|r(S)}$ are equivalent, where
$\sigma(S):s(S)\rightarrow r(S)$, such that $\sigma(S)(x)=r(Sx)$, is the map
induced by $S$. Moreover $d\sigma(S)_*\mu(y)=D_\mu^{-1}(yS)d\mu(y)$.

The {\it Radon-Nikodym problem} for $D\in Z^1(G,{\bf R}_+^*)$, where $G$ is a
locally compact groupoid with Haar system $\lambda$
and compact unit space
$X$, is to determine the set $S_{G,D}(X)$ (or $S_D(X)$ if there is no ambiguity on
$G$) of probability measures
$\mu$ on
$X$ which are quasi-invariant and admit $D$ as their Radon-Nikodym derivative. We
denote by $S(X)$ the set of probability measures on $X$.

\subsection{Cocycles on proper equivalence relations}

Let us first look at cocycles on proper equivalence relations.

\begin{prop}\label{cocycles:proper} Let $R$ be a proper equivalence relation on $X$
defined by $\pi:X\rightarrow\Omega$ and let $c\in Z^1(R,A)$, where $A={\bf R}$ or
${\bf R}_+^*$. Then
\begin{enumerate}
\item $c$ is a coboundary.
\item If $b$ and $b'$ are two potentials for $c$, they differ by a function of the
form $h\circ\pi$, where
$h:\Omega\rightarrow A$ is continuous.
\item If $X$ is compact and $D\in Z^1(R,{\bf R}_+^*)$, there is a unique potential
$\rho\in C(X,{\bf R}_+^*)$ such that $\sum_\omega \rho(x)=1$ for all
$\omega\in\Omega$.
\end{enumerate}
\end{prop}

\begin{proof} For $(i)$ and $(ii)$, we assume that $A={\bf R}$. To prove the first
assertion, we choose a locally finite open cover
$\{V_j\}$ of
$\Omega$, continuous sections $\sigma_j:V_j\rightarrow X$ of $\pi$ and a
partition of unity $\{h_j\}$ subordinate to $\{V_j\}$. We define
$b_j:\pi^{-1}(V_j)\rightarrow{\bf R}$ by $b_j(x)=c(x,\sigma_j(\pi(x))$ and 
$b:X\rightarrow{\bf R}$ by $b(x)=\sum_j h_j(\pi(x))b_j(x)$. Then
$c(x,y)=b(x)-b(y)$.  For the second assertion, we notice that if
$b$ and $b'$ satisfy $b'(x)-b'(y)=b(x)-b(y)$ for all $(x,y)\in R$, then $b'-b$ is
constant on the equivalence classes, hence of the form $h\circ\pi$. For $(iii)$,
we pick an arbitrary potential $\rho'\in C(X,{\bf R}_+^*)$ and define $Z(\omega)=
\sum_\omega \rho'(x)$. Then $Z\in C(\Omega,{\bf R}_+^*)$ and
$\rho=\rho'/Z\circ\pi$ is the desired potential. Because of $(ii)$, two potentials
satisfying this normalization agree.
\end{proof}

Therefore, via the introduction of this normalized potential, we relate a cocycle
$D$ with a family of probability measures along the fibers of the map
$\pi:X\rightarrow\Omega$ or equivalently, a conditional expectation:

\begin{defn}\label{normalized potential} Given a proper equivalence relation $R$
on a compact space $X$ and $D\in Z^1(R,{\bf R}_+^*)$, the potential $\rho=\rho_D\in
C(X,{\bf R}_+^*)$ satisfying the normalization $\sum_\omega \rho(x)=1$ for all
$\omega\in\Omega$ is called the {\it normalized potential} of $D$. It defines a
Markovian (\defnref{positive}) operator $E=E_D: C(X)\rightarrow C(\Omega)$
according to
$E(f)(\omega)=\sum_\omega \rho(x)f(x)$, called the {\it expectation} relative to
$D$.
\end{defn}

\begin{rem} 
\begin{itemize}
\item Our definition carries an abuse of language. It is not $E$ itself but
$\pi^*\circ E$ which is a conditional expectation from $C(X)$ to $C(X)$.
\item I owe to R.~Exel the following observation. The above expectation
$E$ has finite index. Conversely, given a continuous surjection
$\pi:X\rightarrow\Omega$, there exists a conditional expectation
$E:C(X)\rightarrow C(\Omega)$ of finite index if and only if $\pi$ is a local
homeomorphism. Then, there exists $D\in Z^1(R,{\bf R}_+^*)$ such that $E=E_D$.
\end{itemize}
\end{rem}

The Radon-Nikodym problem is easily solved for proper equivalence relations. The
following proposition says that the solutions are the measures on $X$ admitting
$E_D$ as conditional expectation.

\begin{prop} Let $R$ be a proper equivalence relation
on a compact space $X$ and $D\in Z^1(R,{\bf R}_+^*)$. Then
$$S_D(X)=\{\Lambda\circ E_D, \Lambda\in S(\Omega)\},$$
where $E_D$ is the expectation relative to $D$.
More precisely, $E_D^*:C(\Omega)^*\rightarrow C(X)^*$ identifies $S_D(X)$ and
$S(\Omega)$ as compact convex sets.

\begin{proof} Let us first show that, for $\Lambda\in S(\Omega)$,
$\mu=\Lambda\circ E\in S_D(X)$. We denote by $\alpha$ the system of counting
measures on the fibers of $\pi:X\rightarrow\Omega$, by $\lambda_r$ the system of counting
measures on the fibers of $r:R\rightarrow X$ (the first projection) and by
$\lambda_s$ the system of counting
measures on the fibers of $s:R\rightarrow X$ (the second projection). Then, we have
for $f\in C_c(R)$,
$$\mu\circ\lambda_r(f)=\Lambda\circ\alpha(\rho\lambda_r(f))
=\Lambda\circ\alpha\circ\lambda_r((\rho\circ r)f),$$
$$\mu\circ\lambda_s(f)=\Lambda\circ\alpha(\rho\lambda_s(f))
=\Lambda\circ\alpha\circ\lambda_s((\rho\circ s)f).$$
The obvious equality $\alpha\circ\lambda_r=\alpha\circ\lambda_s$ gives the result.
Conversely, let us suppose that $\mu\in S_D(X)$. Then, for every $f\in C(X\times
X)$, we have
$$\int\sum_{\pi(y)=\pi(x)}f(x,y)d\mu(x)=\int\sum_{\pi(x)=\pi(y)}f(x,y)D(x,y)d\mu(y).$$
If we choose $f(x,y)=\rho(y)f(y)$, where $f\in C(X)$, we get
$$\int\sum_{\pi(y)=\pi(x)}\rho(y)f(y)d\mu(x)=\int f(y)d\mu(y),$$
which says that $\pi_*\mu\circ E=\mu$.
\end{proof}

\end{prop}

\begin{lem}\label{decomposition} Let $X\xrightarrow{\pi_1}
X_1\xrightarrow{\pi_{2,1}}X_2$ be surjective local homeomorphisms.
Let $R_1$ [resp.$R_2$] be the equivalence relation defined by $\pi_1$
[resp.$\pi_2=\pi_{2,1}\circ\pi_1$]. Let $D_2\in Z^1(R_2,{\bf R}_+^*)$ with
potential $\rho_2\in C(X,{\bf R}_+^*)$ and $D_1=D_{|R_1}$ with potential $\rho_1$.
Then, there exists a unique
$\rho_{2,1}\in C(X_1,{\bf R}_+^*)$ such that
$\rho_2=\rho_1(\rho_{2,1}\circ\pi_1)$. Moreover, if $\rho_2$ and $\rho_1$ are
normalized, then so is $\rho_{2,1}$, i.e.
$\sum_{\pi_{2,1}(x_1)=x_2}\rho_{2,1}(x_1)=1$ for all $x_2\in X_2$. 
\end{lem}

\begin{proof} The uniqueness of $\rho_{2,1}$ results from the surjectivity of
$\pi_{2,1}$. Since
$\rho_2$ is also a potential for $D_1$, its existence results from
\propref{cocycles:proper} $(ii)$. Let us assume that $\rho_2$ and $\rho_1$ are
normalized. Then for $x_2\in
X_2$, we have
\[\begin{array}{ccc}
1&=&\ \sum_{\pi_2(x)=x_2}\rho_2(x)\\
&=&\ \sum_{\pi_{2,1}(x_1)=x_2}\sum_{\pi_1(x)=x_1}\rho_1(x)\rho_{2,1}\circ\pi_1(x)\\
&=&\ \sum_{\pi_{2,1}(x_1)=x_2}\rho_{2,1}(x_1)\sum_{\pi_1(x)=x_1}\rho_1(x)\\
&=&\ \sum_{\pi_{2,1}(x_1)=x_2}\rho_{2,1}(x_1)\\
\end{array}\]
\end{proof}

\subsection{Cocycles on AP equivalence relation}

Let us consider now the case of an AP equivalence relation $R$ on a compact space
$X$ with a defining sequence $(X_n)$. We use the notation of \defnref{AP}; in
particular, we denote by
$R_n$ the proper equivalence relation defined by $\pi_n:X\rightarrow X_n$. Given
$D\in Z^1(R,{\bf R}_+^*)$, we consider the sequence of its restrictions
$D_n=D_{|R_n}\in Z^1(R_n,{\bf R}_+^*)$. For each $n$, we pick a potential
$\rho_n$ for $D_n$. From
\lemref{decomposition}, we obtain
$\rho_{n,n-1}\in C(X_{n-1},{\bf R}_+^*)$ such that
$\rho_n=\rho_{n-1}(\rho_{n,n-1}\circ\pi_{n-1})$. Equivalently,
$\rho_n(x)/\rho_{n-1}(x)$ depends only on $\pi_{n-1}(x)$. In particular, we can
choose the sequence of normalized potentials; in that case $\rho_{n,n-1}$ is also
normalized.

\begin{defn} Let $R$ be an AP equivalence relation on a compact space
$X$ defined by a sequence $(X_n)$.
\begin{enumerate}
\item A {\it compatible sequence of potentials} is a sequence
$(\rho_n\in C(X,{\bf R}_+^*))$ such that
$\rho_n(x)/\rho_{n-1}(x)$ depends only on $\pi_{n-1}(x)$. Equivalently, a
compatible sequence of potentials is defined by an initial potential 
$\rho_0\in C(X,{\bf R}_+^*)$ and a sequence of
local potentials
$\rho_{n,n-1}\in C(X_{n-1},{\bf R}_+^*)$; then $\rho_n$ is given by
$$\rho_n=(\rho_0)(\rho_{1,0})(\rho_{2,1}\circ\pi_1)\ldots(\rho_{n,n-1}\circ\pi_{n-1}).$$
If for all
$n$ and  for all $x_n\in X_n$, $\sum_{x_n}\rho_n(x)=1$, we speak of a sequence of
normalized potentials.
\item Given $D\in Z^1(R,{\bf R}_+^*)$, a sequence
of potentials
$\rho_n$ of $D_n=D_{|R_n}$ is called a sequence of {\it potentials} of
$D$. The associated sequence of funtions
$(\rho_{n,n-1})$ is called a sequence of {\it local potentials} of $D$. In
particular, we can consider the sequence of {\it normalized potentials} of $D$.
\end{enumerate}
\end{defn}

Note that a compatible sequence of potentials $(\rho_n)$ defines a cocycle $D\in
Z^1(R,{\bf R}_+^*)$ such that $D(x,y)=\rho_n(x)/\rho_n(y)$
for $(x,y)\in R_n$ and that every cocycle $D\in
Z^1(R,{\bf R}_+^*)$ is defined by a compatible sequence of potentials.

\begin{rem} The sequence of normalized potentials of a cocycle $D\in Z^1(R,{\bf
R}_+^*)$ is unique. However, it is sometimes advantageous to consider arbitrary
sequences of potentials. For example, a change of initial potential amounts to
replacing the cocycle by a cohomologous cocycle. Although we restrict the
discussion below to the normalized potentials, the theory of dimension groups and
their state space applies as well to the general case.
\end{rem}

\begin{prop} The correspondence wich associates to a cocycle $D$ its sequence of
normalized potentials $(\rho_n)$ [resp. its sequence of normalized local
potentials $(\rho_{n,n-1})$] is a bijection from $Z^1(R,{\bf R}_+^*)$ onto
the set of compatible sequences of normalized potentials [resp. onto the set of
sequences of normalized local potentials].
\end{prop}

The normalized potentials $\rho_n$ define
conditional expectations
$E_n:C(X)\rightarrow C(X_n)$ and the normalized local potentials $(\rho_{n,n-1})$
define conditional expectations
$$E_{n,n-1}: C(X_{n-1})\rightarrow C(X_n).$$ 
For $m\le n$, we set:
$$E_{n,m}=E_{n,n-1}E_{n-1,n-2}\ldots E_{m+1,m}.$$
Then, we have $E_n=E_{n,0}$.

\begin{prop}\label{martingale ppty} Let $R=\cup R_n$ be an AP equivalence relation
on a compact space $X$, $D\in Z^1(R,{\bf R}_+^*)$ and $D_n=D_{|R_n}\in
Z^1(R_n,{\bf R}_+^*)$ as above. Then, the sequence of conditional expectations
$(P_n=\pi_n^*E_n)$ is a (reversed) martingale, i.e. it satisfies
$P_mP_n=P_nP_m=P_n$ for $m\le n$.
\end{prop}

\begin{proof} This is an immediate consequence from the fact that
$E_m\pi_m^*$ is the identity map of
$C(X_m)$ and of the factoriations $\pi_n=\pi_{n,m}\circ\pi_m$ and $E_n=E_{n,m}E_m$.
\end{proof}

Thus we obtain a sequence of Markovian operators
$$\underline E=(E_{n,n-1}:C(X_{n-1})\rightarrow C(X_n))$$
and we can apply the
results of the Appendix. We denote by ${\cal E}={\cal E}(\underline E)$ its
dimension group and by $S=S(\underline E)$ its state space. Recall that it is a
non-empty compact convex Choquet simplex. We denote by
$E_{\infty,n}:C(X_n)\rightarrow {\cal E}$ the canonical morphisms; they are
surjective. In particular, we write $E_\infty=E_{\infty,0}:C(X)\rightarrow {\cal
E}$. We have a convenient
description of the state space $S(\underline E)$. Indeed, given a reversed
martingale $(P_n)$, Brown and Dooley define in \cite{bd:g-measures}  a $G$-measure
as a probability measure $\mu$ on
$X$ satisfying
$\mu= (P_n)^*\mu$ for all $n$. On the other hand, by
\defnref{state space}, S(\underline E) is the set of sequences
$(\mu_n)$, where
$\mu_n\in S(X_n)$ and $\mu_{n-1}=\mu_n\circ E_{n,n-1}$. For such a sequence,
$\mu=\mu_0=E_n^*\mu_n$ does not depend on $n$ and is a $G$-measure. Conversely, one
recovers the sequence $(\mu_n)$ from the $G$-measure $\mu$ by setting
$\mu_n=(\pi_n)_*\mu$. Thus:

\begin{lem} The map $(\mu_n)\mapsto \mu_0$ identifies the state space
$S(\underline E)$ to the set of $G$-measures.
\end{lem} 

Moreover, we have seen previously that the measures of the form $E_n^*\mu_n$, where
$\mu_n\in S(X_n)$ are exactly the probability measures admitting $D_n$ as
Radon-Nikodym derivative. This gives the next proposition. In the context of
equilibrium states in statistical mechanics, the 
equivalence of the definitions of a Gibbs measure as a $G$-measure or as a
quasi-invariant measure is well known (\cite{kel:equilibrium}, Theorem
5.2.4.$(a)$).

\begin{prop}\label{state space} Let $(X,R)$, $D\in
Z^1(R,{\bf R}_+^*)$, $D_n$, $E_n$ and $\underline E$ as above.
\begin{enumerate}
\item The sequence
$(S_{D_n}(X))$ is decreasing and its intersection is $S_D(X)$. 
\item  The isomorphisms $E_n^*:S(X_n)\rightarrow S_{D_n}(X)$ induce an isomorphism
of compact convex sets from $S=S(\underline E)$ onto $S_D(X)$. In other words, the
quasi-invariant probability measures admitting $D$ as Radon-Nikodym derivative are
exactly the $G$-measures of the martingale $(P_n=\pi_n^*E_n)$.
\end{enumerate}
\end{prop}

\begin{proof} $(i)$ A measure $\mu\in S(X)$ belongs to $S_{D_n}(X)$ if and only if
it is quasi-invariant with respect to $R_n$  and admits $D_n$ as its
Radon-Nikodym derivative. Since $ R_{n-1}\subset R_n$ and
$D_{n-1}={D_n}_{|R_{n-1}}$, a measure $\mu\in S(X)$ which belongs to $S_{D_n}(X)$
also belongs to $S_{D_{n-1}}(X)$. Moreover $\mu\in S(X)$ belongs to $S_D(X)$ if
and only if, for all $n$, it is quasi-invariant with respect to $R_n$ with
Radon-Nikodym derivative
$D_n$. Therefore, the intersection of the $(S_{D_n}(X))$'s is $S_D(X)$.

$(ii)$ Let $(\mu_n)\in S(\underline E)$ and let $\mu=E_n^*\mu_n$ the associated
$G$-measure. Then $\mu$ belongs to
$S_{D_n}(X)$ for all
$n$ and therefore belongs to $S_D(X)$. Conversely, if $\mu$ belongs to
$S_D(X)$, we define $\mu_n={\pi_n}_*\mu\in S(X_n)$ for all $n$. Then the sequence
$(\mu_n)$ belongs to $S(\underline E)$ and $\mu=E_n^*\mu_n$. Therefore, the map
which sends
$(\mu_n)\in S(\underline E)$ into $\mu_0\in S(X)$ is an isomorphism of
$S(\underline E)$ onto $S_D(X)$.
\end{proof}

We reformulate \lemref{key lemma} of the Appendix in this
setting.

\begin{cor}\label{martingale convergence} Let $(X,R)$, $D\in
Z^1(R,{\bf R}_+^*)$, $E_n$ be as above. For $f\in C(X)$, the following
conditions are equivalent:
\begin{enumerate}
\item  $\|E_n(f)\|$ tends to $0$ as $n$ tends to $\infty$,
\item  $\mu(f)=0$ for all $\mu\in S_D(X)$.
\end{enumerate}
\end{cor}

With the notation of the Appendix, we denote by $Aff(S)$ the space of continuous
functions on the compact convex space $S$ and by
$\theta:{\cal E}\rightarrow Aff(S_D(X))$ the evaluation map. Then
$\theta_0=\theta\circ E_\infty:C(X)\rightarrow Aff(S_D(X))$ is the evaluation map
$\theta_0(f)(\mu)=\mu(f)$ for $\mu\in S_D(X)$.

\begin{prop}\label{inf} Let $(X,R)$ be an AP equivalence relation and let $D$ be a
continuous cocycle in
$Z^1(R,{\bf R}_+^*)$. Then the evaluation map $\theta_0$ induces an isomorphism of
ordered Banach spaces with order-unit from $C(X)/Ker\theta_0$ onto $Aff(S_D(X))$. 
\end{prop}

\begin{proof} This results from \propref{affine representation}.
\end{proof}

In the next lemma, $(\pi_\infty)_*$ is the restriction map given by the
inclusion $C(X)^R\subset C(X)$.

\begin{lem}\label{direct sum}  Let $(X,R)$ be an AP equivalence relation and let
$D$ be a continuous cocycle in
$Z^1(R,{\bf R}_+^*)$. Then,
\begin{enumerate}
\item $C(X)^R\cap Ker\theta_0=\{0\}$ and
\item $(\pi_\infty)_*:S_D(X)\rightarrow S(X_\infty)$ is surjective.
\end{enumerate}
\end{lem}

\begin{proof}
$(i)$ Let $f$ be a non-zero element of $C(X)^R$. There exists
$x\in X$ be such that $f(x)=a\not=0$. The sequence of probability measures
$(E_n^*\delta_{\pi_n(x)})$ admits at least one limit point $\mu$ is $S(X)$. Since
$E_n^*\delta_{\pi_n(x)}$ belongs to $S_{D_n}(X)$, $\mu$ belongs to $S_D(X)$. Since
$$E_n^*\delta_{\pi_n(x)}(f)=E_n(f)(\pi_n(x))=f(x)=a$$
for all $n$, $\mu(f)=a\not=0$.

$(ii)$ results from $(i)$ and the Hahn-Banach theorem (one can use for example
Proposition 4.2 of \cite{goo:poag}).
\end{proof}

\begin{prop}\label{equicontinuous} Let $(X,R=\cup R_n)$ and $D\in Z^1(R,{\bf
R}_+^*)$ be as above. The following assertions are equivalent:
\begin{enumerate}
\item $C(X)=C(X)^R+Ker\theta_0$.
\item The restriction of $\theta_0$ to $C(X)^R$ is an isomorphism from $C(X)^R$
onto
$Aff(S_D(X))$.
\item The restriction map $(\pi_\infty)_*$ is an isomorphism
from
$S_D(X)$ onto $S(X_\infty)$.
\item There is a conditional expectation $E_\infty:C(X)\rightarrow C(X_\infty)$
such that
$$S_D(X)=\{\Lambda\circ E_\infty: \Lambda\in S(X_\infty)\}.$$
\item For all $f\in C(X)$, $(\pi_n^*E_n(f))$ converges uniformly.
\item For all $f\in C(X)$, $(\pi_n^*E_n(f))$ is equicontinuous.
\end{enumerate}
\end{prop}

\begin{proof} $(i)\Leftrightarrow (ii)$ This results from
\propref{inf} and
\lemref{direct sum}.

$(ii)\Leftrightarrow (iii)$ We can identify the state space of $Aff(S_D(X))$ with
$S_D(X)$ and the state space of $C(X)^R$ with $S(X_\infty)$. The map from
$S_D(X)$ into $S(X_\infty)$ induced by $\theta_0:C(X)^R\rightarrow Aff(S_D(X))$
is the restriction map
$(\pi_\infty)_*$. Therefore, if $(ii)$ holds, this restriction map is an
isomorphism of convex compact sets. Conversely,
let us assume that the restriction map $(\pi_\infty)_*:S_D(X)\rightarrow
S(X_\infty)$ is an isomorphism of compact convex sets.
It induces an isomorphism of affine spaces from
$Aff(S(X_\infty))=C(X^\infty)$ onto $Aff(S_D(X))$. But this map coincides with
$\theta_0$.

$(i)\Rightarrow (iv)$ We let $P$ be the projection of $C(X)$ onto $C(X)^R$ along
$Ker\theta_0$ and $E_\infty:C(X)\rightarrow C(X_\infty)$ be the composition of $P$
and the isomorphism $C(X^R)\rightarrow C(X_\infty)$. Then, for all $\mu\in S_D(X)$
and $f\in C(X)$, we have the equalities
$$<\mu,\pi_\infty^*E_\infty(f)>=<\mu, Pf>=<\mu,f>.$$
and therefore $\mu=(\pi_\infty)_*\circ E_\infty$. Since the restriction map
$(\pi_\infty)_*$ is a bijection from
$S_D(X)$ onto
$S(X_\infty)$, we have the result
$(iv)$.

$(iv)\Rightarrow (v)$ Let $f\in C(X)$. Then, using the characterization of $S_D(X)$
given in $(iv)$, we see that
$f-\pi_\infty^*E_\infty(f)$ belongs to $Ker\theta_0$. According to 
\corref{martingale convergence},
this implies that
$\|\pi_n^*E_n(f)-\pi_\infty^*E_\infty(f))\|=\|E_n(f-\pi_\infty^*E_\infty(f))\|$
tends to
$0$.

$(v)\Rightarrow (vi)$ is clear.

$(vi)\Rightarrow (i)$ To ease the notation, we introduce
$P_n=\pi_n^*E_n:C(X)\rightarrow C(X)$. Let
$f\in C(X)$. Since
$(P_n(f))$ is equicontinuous and bounded, there is a  subsequence
$(P_{n_k}(f))$ converging to $g\in C(X)$. For a fixed $n$,
$(P_nP_{n_k}(f))$ converges to $P_n(g)\in C(X)$ by continuity of $P_n$. But since
$P_nP_{n_k}(f)=P_{n_k}(f)$ for $k$ sufficiently large, the sequence converges
also to $g$, hence $g=P_n(g)$ and belongs to $\pi_n^*(C(X_n)$. Thus $g$ belongs
to $C(X)^R$. It remains to show that $f-g$ belongs to $Ker\theta_0$. According to
\corref{martingale convergence}, it suffices to check that $\|E_n(f-g)\|$ tends to
$0$. But this is clear, since for $n\ge n_k$, we have the inequality 
$$\|E_n(f-g)\|\le \|P_{n_k}(f)-g\|.$$
\end{proof}

Let us specialize the above proposition to the irreducible case, i.e. when
$X_\infty$ is reduced to one point. It will give a
necessary and sufficient condition for the uniqueness of the solution of the
Radon-Nikodym problem.

The following corollary of \propref{martingale convergence}, which is well known
(see for example \cite{bd:g-measures} or Proposition 1 of \cite {fan:Ruelle}),
will be our basic tool to show unique ergodicity.

\begin{cor}\label{unique ergodicity} Let
$(X,R=\cup R_n)$ and
$D\in Z^1(R,{\bf R}_+^*)$ be as above. The following assertions are equivalent.
\begin{enumerate}
\item $D$ is {\it uniquely ergodic}, i.e.  $S_D(X)$ contains one element.
\item $(X,R)$ is irreducible and $D$ is equicontinuous (i.e. for all $f\in C(X)$,
$(\pi_n^*E_n(f))$ is equicontinuous).
\item For all $f\in C(X)$, $(\pi_n^*E_n(f))$ converges uniformly to a constant
function.
\item For all $f\in C(X)$, the variation $var_+(E_n(f))$ of $E_n(f)$ (i.e. the
difference between its maximum and its minimum) tends to
$0$.
\end{enumerate}
\end{cor}

\begin{proof} We first observe that $(i)$ or $(ii)$ imply that
$C(X)^R={\bf R}1$. For  $(i)$, this is a consequence of
\lemref{direct sum}. For $(ii)$, this comes from the fact that $\pi_n^*E_n$ acts
as the identity on $C(X)^R$. Then, one can see that these first three conditions
are a reformulation of the conditions of
\propref{equicontinuous} under the assumption that $R$ is irreducible. The
equivalence of $(i)$ and $(iv)$ is a particular case of \corref{CNS unique state}.
\end{proof}

\section{Quasi-product cocycles on AF relations.}  Let
$(V,E)$ be a Bratteli diagram. We denote by $X=X(V,E)$ its infinite path space, by
$X_n$ the space of infinite paths starting at level $n$ and
by $R=R(V,E)=\cup R_n$ the tail equivalence relation on $X$. A function
$\Phi:E\rightarrow {\bf R}_+^*$ defines a cocycle $D\in Z^1(R,{\bf R}_+^*)$
according to the formula
$$D(x,y)=\lim_{n\to\infty} {\Phi(x_1)\Phi(x_2)\ldots
\Phi(x_n)\over\Phi(y_1)\Phi(y_2)\ldots
\Phi(y_n)}.$$
(Note that the sequence is stationary.)

\begin{defn} Given a Bratteli diagram $(V,E)$, a cocycle $D$ on the tail
equivalence relation $R(V,E)$ of its infinite path space $X(V,E)$ is called a
{\it quasi-product cocycle} if it is of the above form.
\end{defn}

\begin{rem}
It is shown in \cite{ren:AF} that every continuous cocycle on an AP equivalence
relation on a compact totally disconnected space is cohomologous to a quasi-product
cocycle relative to some Bratteli diagram. However, this result gives little
information about the Bratteli diagram, in particular about its growth. We shall
characterize later the cocycles which are cohomologous to a quasi-product
cocycle relative to some contraction of a given Bratteli diagram.
\end{rem}
We assume from now on that the set of vertices $V(n)$ of each level $n$ is
finite (and as before, that every vertex
emits  finitely many, but at least one, edges and that every vertex of a
level $n\ge 1$ receives at least one edge). Then the space of infinite paths $X$
of the diagram is compact. Let $D$ be a quasi-product cocycle given by
$\Phi:E\rightarrow {\bf R}_+^*$. Let us apply to
$D$ the general approximation procedure. The normalized potential of
$D_n=D_{|R_n}$ is the function
$\rho_n\in C(X)$ given by
$$\rho_n(x)=\Phi(x_1)\Phi(x_2)\ldots \Phi(x_n)/Z_n(r(x_n)) ,$$
where we have introduced the normalization factor
$$Z_n(v)=\sum \Phi(y_1)\Phi(y_2)\ldots \Phi(y_n)$$ where $v\in V(n)$ and the sum is
over all the finite paths $y_1y_2\ldots y_n$ of length $n$ ending at $v$. The
sequence of normalized local potentials of $D$ is given by
$$\rho_{n,n-1}(y_nx_{n+1}\ldots)=Z_n\circ r(y_n)^{-1}\Phi(y_n) Z_{n-1}\circ
s(y_n).$$  
Note that $\rho_{n,n-1}$ depends only
on $y_n$. Let us point out this elementary property in the following proposition.
\begin{prop}\label{quasi-product cocycle} Let $(X,R)$ be the AF equivalence
relation defined by the Bratteli diagram $(V,E)$ and let $D\in Z^1(R, {\bf
R}_+^*)$. Then, the following conditions are equivalent:
\begin{enumerate}
\item $D$ is a quasi-product cocycle relative to $(V,E)$.
\item For all $n$, its normalized local potential
$\rho_{n,n-1}$ depends only on the first variable $x_n$.
\item For all $n$, its normalized potential
$\rho_n$ depends only on the first $n$ variables: $x_1x_2\ldots x_n$.
\item $D$ admits an initial potential $\rho_0\equiv 1$ and a sequence of local
potentials
$(\rho_{n,n-1})$ such that
$\rho_{n,n-1}(x_n\ldots)$ depends only on $x_n\in E(n)$.
\end{enumerate}
\end{prop}
\begin{proof} The equivalence of $(i)$ and $(iv)$ is clear. It is also clear that
$(iv)$ implies $(ii)$, and that $(ii)$ implies $(iii)$, because
$$\rho_n(x)=\rho_{1,0}(x)\ldots\rho_{n,n-1}\circ\pi_{n-1}(x).$$
Finally, $(iii)$ implies $(iv)$ because 
$\rho_{n,n-1}\circ\pi_{n-1}(x)=\rho_n(x)/\rho_{n-1}(x)$.
\end{proof}
The conditional
expectation
$E_n:C(X)\rightarrow C(X_n)$ is given by
$$E_n(f)(x_{n+1}\ldots)=\sum \rho_n(y_1\ldots y_nx_{n+1}\ldots)f(y_1\ldots
y_nx_{n+1}\ldots)$$
where the sum is
over all the finite paths $y_1y_2\ldots y_n$ of length $n$ ending at $s(x_{n+1})$.
The conditional expectation $E_{n,n-1}:C(X_{n-1})\rightarrow C(X_n)$ is given by
$$E_{n,n-1}(f)(x_{n+1}\ldots)=\sum
\rho_{n,n-1}(y_nx_{n+1}\ldots)f(y_nx_{n+1}\ldots)$$ where the sum is over
all the $y_n\in E(n)$ ending at $s(x_{n+1})$. 
Thus we have as usual a sequence
$\underline E=(E_{n,n-1}:C(X_{n-1})\rightarrow C(X_n))$ of Markovian operators. On
the other hand, $\Phi$ defines a sequence of matrices
$\underline A=(A_n:C(V(n-1))\rightarrow C(V(n)))$, where the coefficients of  $A_n$
are $A_n(w,v)=\sum\Phi(e)$, the sum being over all the $e\in E(n)$ starting at
$v\in V(n-1)$ and ending at $w\in V(n)$. The state space $S(\underline A)$ of
this sequence is defined in the Appendix. It is the state space of the inductive
limit dimension group. By definition, it consists of sequences $(\rho_n)$ of
positive numbers such that 
\begin{eqnarray}
\rho_{n-1}(v)&=&\sum_{s(e)=v}\Phi(e)\rho_n\circ r(e) \quad
(n=1,2,\ldots\ v\in V(n-1))\\ 
1&=&\sum_{V(0)}\rho_0(v).
\end{eqnarray}
Such a sequence $(\rho_n)$
defines a {\it Markov measure} $\mu$ such that
$$\mu(Z(x_1\ldots x_n))=\rho_0(s(x_1))p_1(x_1)\ldots p_n(x_n)$$
where $p_n(e)=(\rho_{n-1}\circ s(e))^{-1}\Phi_n(e)\rho_n\circ r(e)$ for
$n=1,2,\ldots, e\in E(n)$. 

\begin{prop}\label{quasi-product cocycle} (\cite{ren:AF}, Proposition 3.3) Let
$D$ be a quasi-product cocycle. The above construction defines an isomorphism from
the state space $S(\underline A)$ of
the sequence of matrices $\underline A=(A_n:C(V(n-1))\rightarrow C(V(n)))$ to
the space
$S_D(X)$ of solutions of the
$(RN)$ equation $D_\mu=D$. In other words, the
quasi-invariant probability measures admitting $D$ as
Radon-Nikodym derivative are exactly the above Markov measures.
\end{prop}

\begin{proof} We refer the reader to \cite{ren:AF}. It easy to check that the
Markov measure $\mu$ constructed above belongs to $S_D(X)$. It is shown in
\cite{ren:AF} that  conversely every $\mu\in S_D(X)$ is of that form. The maps so
defined are continuous and preserve convex combinations.
\end{proof}

Note that the maps
$J_n:C(V(n))\rightarrow C(X_n)$ defined by $J_n(f)(x_{n+1}\ldots )=Z_n\circ
s(x_{n+1})^{-1}f\circ s(x_{n+1})$ yield a morphism from $\underline A$ to
$\underline E$, hence a morphism of their dimension groups $J:{\cal E}(\underline
A)\rightarrow {\cal E}(\underline E)$. This morphism induces the above isomorphism
of their state spaces. However $J$ itself is not necessarily an isomorphism. For
example, let the Bratteli diagram $(V,E)$ be a tree (and $\Phi:E\rightarrow{\bf
R}_+^*$ be identically one). Then
$\pi_{n.n-1}:X_{n-1}\rightarrow X_n$ is a bijection and
$E_{n,n-1}:C(X_{n-1})\rightarrow C(X_n)$ is the transposed map $\pi^*_{n.n-1}$.
The inductive limit $\underline E$ can be identified with $C(X)$. On the other hand
the image of $J$ consists of locally constant functions.

Taking into account this proposition, \corref{CNS unique state} and
\corref{CS unique state} become:

\begin{cor}\label{AF unique state} Let $R$ be the tail equivalence relation on the
infinite path space
$X$ of a Bratteli diagram $(V,E)$ and let $D\in Z^1(R,{\bf R}_+^*)$ be a
quasi-product cocycle defined by $\Phi:E\rightarrow {\bf R}_+^*$. Define the matrix
$A_n(w,v)=\sum\Phi(e)$, the sum being over all the $e\in E(n)$ starting at
$v\in V(n-1)$ and ending at $w\in V(n)$ as
above. 
\begin{enumerate}
\item A necessary and sufficient condition for $S_D(X)$ to have exactly one element
is that, for any fixed $m$ and $v\in V(m)$, the variation of the function $w\in
V(n)\mapsto B_n\ldots B_{m+1}(v,w)$ goes to 0 when $n$ goes to infinity, where
$u_n(w)=\sum_{v\in V(0)} A_n\ldots
A_1(w,v)$ and
$B_n(w,v)=u_n(w)^{-1}A_n(w,v)u_{n-1}(v)$ for $w\in V(n)$ and $v\in V(n-1)$.
\item  A sufficient condition for $S_D(X)$ to have exactly one element
is that the serie $\sum\epsilon_n$ diverges, where $\epsilon_n$ is the ratio
of the smallest element of $A_n$ over its largest element.
\end{enumerate}
\end{cor}

The most studied quasi-product cocycle is $D\equiv 1$. Then, the elements of
$S_D(X)$ are the invariant (with respect to the tail equivalence relation)
probability measures on the path space $X$ of the Bratteli diagram $(V,E)$ (they
are called central measures in \cite{vk:asymptotic}). They correspond to the
tracial states of the AF algebra of the Bratteli diagram. We choose
$\Phi\equiv 1$. The matrices
$A_n$ are the adjacency matrices of the graph. They have coefficients in
$\bf N$ and the inductive limit of the $(A_n:C(V(n-1,{\bf Z}))\rightarrow
C(V(n),{\bf Z}))$'s is the usual dimension group of the Bratteli diagram. Then
\propref{quasi-product cocycle} gives the well known correspondence
between invariant probability measures and states of the dimension group. In that
case, parts of \corref{AF unique state} appear in the work
\cite{toe:AF} of A.~T\"or\"ok. The necessary and
sufficient condition
$(i)$ also appears in \cite{vk:asymptotic}. Here are a few examples.

\begin{ex}{\it An example of Fack and Mar\'echal}. In their work \cite{fam:sym} on
the symmetries of UHF algebras, T.~Fack and O.~Mar\'echal study the dimension group
defined by the sequence of matrices $\underline A=(A_n)$, where 
$$\displaystyle\,
A_n=\begin{pmatrix} p_n & r_n\\
r_n & p_n
\end{pmatrix}$$
and $(p_n)$ and $(r_n)$ are two sequences of integers such that $0<r_n\le p_n$
for all $n$. Let $\epsilon_n=r_n/p_n$. According
to the above, a sufficient condition for ${\cal E}(\underline A)$ having a 
unique state is that
$\,\sum\epsilon_n=\infty\,$. Fack and Mar\'echal show by an expli\-cit computation
of the dimension group that this condition is necessary. This can also be deduced
from \corref{AF unique state} $(i)$.

\end{ex}

\begin{ex}{\it Pascal's triangle}.
The simple random walk on $\bf Z$ provides another example. We let
$X=\prod_1^\infty\{0,1\}$ be the space of increments. We introduce
$X_n=\{0,1,\ldots,n\}\times\prod_{n+1}^\infty\{0,1\}$ and
$\pi_n:X\rightarrow
X_n$ defined by
$$\pi_n(x_1x_2\ldots)=(x_1+\ldots+x_n,x_{n+1}x_{n+2}\ldots).$$ Its first
coordinate is the position of the walker at time $n$ (assuming that his initial
position is 0). We let $R$ be the AP equivalence relation on $X$ defined by the
$\pi_n$'s. Note that
$(X,R)$ admits the infinite Pascal triangle as Bratteli diagram and that
$D\equiv 1$ is the  quasi-product cocycle defined by the function
$\phi\equiv 1$. It is known that
the state space $S_1$ of invariant probability measures is isomorphic to the
space of probability measures on $[0.1]$; see for example Section VII.4 of
\cite{fel:proba}, \cite {ker:young} or the Appendix of
\cite{ren:approach}, which contains an explicit computation of the dimension group
of the infinite Pascal triangle. The
following comment of the proof is inspired by Section 5 of \cite{wasa:thesis}.
First, it is immediate to check that for
$t\in [0,1]$, the product measure
$$\mu_t=\prod_1^\infty((1-t)\delta_0+t\delta_1)$$
is invariant and that $\mu_t(Z(n,k))=C_n^k t^k(1-t)^{n-k}$, where $Z(n,k)$ is the
set of paths having position $k$ at time $n$. An elementary estimate using the
expansion of a polynomial of degree not greater than
$n$ in the basis of Bernstein polynomials $\{t^k(1-t)^{n-k}, k=0,1,\ldots,n\}$
shows that for all $f\in C(X)$, 
$$\lim_{n\to\infty}\sup_{(k,x)\in X_n} |E_n(f)(k,x)-\hat f(k/n)|=0,$$
 where
$E_n(f)(k,x)$ is the average of the $f(a_1\ldots a_n x)$ over all $a_1\ldots a_n$
such that $a_1+\ldots +a_n=k$ and $\hat f(t)=\mu_t(f)$. Using \corref{martingale
convergence}, one deduces that $\hat f=0$ if and only if $\mu(f)=0$ for all
invariant probability measures $\mu$. Therefore, the measures $\mu_t, t\in[0,1]$
are exactly the extremal elements of $S_1$. Since the equivalence relation $R$ is
irreducible (it has dense orbits), $C(X)^R={\bf C}1$ and the condition $(ii)$ of
\propref{equicontinuous} is not realized. One can also see that for
$x\in X$,
$\pi_n^*E_n(f)(x)$
converges to $\hat f (t)$ for all $f\in C(X)$  if and only if
the path
$x=x_1x_2\ldots$ has the property that
$(x_1+\ldots+x_n)/n$ tends to $t$. The strong law of large numbers (or also the
martingale convergence theorem) says that, with respect to the measure $\mu_t$,
this set is conull.
\end{ex}

\begin{ex}{\it Stationary cocycles}. The Bratteli diagram $(V,E)$ is
called stationary if for all $n\ge 1$, $V(n)=V(0)$ and $E(n)=E(1)$. In that
case, the one-sided shift $T(x_1x_2\ldots)=x_2x_3\ldots$ acts on the infinite
path space $X$ of the diagram. Then, $\Phi:E\rightarrow{\bf R}_+^*$ is called
stationary if it does not depend on the level $n$. The associated quasi-product
cocycle is also called stationary. The sufficient condition $(ii)$ of \corref{AF
unique state} is always satisfied for a stationary quasi-product cocycle $D$.
Therefore $S_D(X)$ is reduced to one element.
\end{ex}

\section{Beyond quasi-product cocycles.}

As mentioned previously, every continuous cocycle $D\in Z^1(R,{\bf R}_+^*)$ of an
AF equivalence relation $R$ on a compact space $X$ is cohomologous to a
quasi-product cocycle with respect to some Bratteli diagram. However, in general we
do not have sufficient information on the Bratteli diagram in order
to apply \corref{AF unique state}. When the cocycle
$D$ satisfies an appropriate condition of equicontinuity relative to a given Bratteli
diagram, it is approximately a quasi-product cocycle with respect to the same Bratteli
diagram (or a contraction of it); then one can use the technique of
\corref{AF unique state} to obtain unique ergo\-dicity. Moreover, the
idea of approximating a cocycle by a quasi-product cocycle is also fruitful for
arbitrary AP equivalence relations. The definitions below are adapted from
\cite{ren:AF}.

\begin{defn}\label{picover} Let $\pi:X\rightarrow Y$ be a surjective local
homeomorphism.  We say that an open subset $V$ of $Y$ is {\it well-covered} if
$\pi^{-1}(V)$ is the disjoint union of a family of open sets
$\{U_i, i\in I\}$ which all map homeomorphically onto $V$. By
definition, a $\pi$-{\it cover} will consist of a cover ${\cal V}$ of $Y$ by
well-covered open subsets and for each $V\in {\cal V}$ a partition of $\pi^{-1}(V)$
by open subsets $U$ of
$X$ which all map homeomorphically onto $V$. We denote by
${\cal U}$ the cover of $X$ by these open sets $U$.
\end{defn} 
When $X$ is
compact,
$\pi$ admits finite $\pi$-covers.
Moreover, if an open cover $\cal W$ of $X$ is given, we can construct our
$\pi$-cover such that $\cal U$ strictly refines $\cal W$ (notation:
${\cal U}\prec{\cal W}$), in
the sense that for each
$U\in {\cal U}$, there exists $W\in {\cal W}$ such that $\overline
U\subset W$.


In the sequel, we shall only consider finite
$\pi$-covers. Recall from \cite{bou:topo} that a finite open cover ${\cal U}$ of a
compact space $X$ defines the entourage
$$\Delta_{\cal U}=\cup_{\cal U}U\times U$$
of the canonical uniform structure of $X$.

Given $\Delta\subset X\times X$ and $\varphi\in C(X,{\bf R})$, we define the
additive variation of $\varphi$ over $\Delta$ as 
$$var_+(\varphi,\Delta)=\sup_{\Delta}|{\varphi(x)-\varphi(x')}|.$$
Similarly, given  $g\in C(X,{\bf R}_+^*)$, we define the
multiplicative variation of $g$ over $\Delta$ as 
$$var_*(g,\Delta)=\sup_{\Delta}|{g(x)\over g(x')}-1|.$$

\begin{lem}\label{normalization} Let $\pi:X\rightarrow Y$ be a surjective local
homeomorphism. Let
$({\cal V}, {\cal U})$ be a
$\pi$-cover with associated entourages
$\delta=\Delta_{\cal V},\Delta=\Delta_{\cal U}$. Consider a potential $g\in C(X,{\bf
R}_+^*)$ and its normalized potential $\rho=(Z\circ\pi)^{-1} g$, where
$Z(y)=\sum_{\pi^{-1}(y)} g(x)$ for $y\in Y$. If $var_*(g,\Delta)<\epsilon<1$, then
$var_*(\rho,\Delta)<2\epsilon(1-\epsilon)^{-1}$.
\end{lem}

\begin{proof} For all $(x,x')\in
\Delta$, we have
$$(1-\epsilon) g(x')< \rho(x) < (1+\epsilon) g(x').$$ Let
$(y,y')\in\delta$. We have by construction a bijection
$x\in
\pi^{-1}(y)\mapsto x'\in \pi^{-1}(y')$ such that
$(x,x')\in \Delta$. Summing above inequalities above $\pi^{-1}(y)$,
we obtain $$(1-\epsilon)Z(y')< Z(y) < (1+\epsilon)
Z(y').$$
Therefore, for all $(x,x')\in
\Delta$, we have
 $${1-\epsilon\over 1+\epsilon}\, \rho(x')<  \rho(x) <
{1+\epsilon\over 1-\epsilon}\, \rho(x')$$ 
and the above inequality.
\end{proof}

We keep the same notation as above. We introduce the positive linear map
$E: C(X)\rightarrow C(Y)$ such that 
$$E(f)(y)=\sum_{\pi^{-1}(y)} \rho(x)f(x).$$

\begin{lem}\label{Walters} (cf. \cite{wal:78}, Lemma 1) Let $\pi:X\rightarrow
Y$ be a surjective local homeomorphism. Let
$({\cal V}, {\cal U})$ be a
$\pi$-cover with associated entourages
$\delta,\Delta$. Consider a normalized potential $\rho\in C(X,{\bf R}_+^*)$
and its expectation $E: C(X)\rightarrow C(Y)$. Then, for every
$f$ in $C(X)$, 
$$var_+(E(f),\delta)\le var_+(f,\Delta)+\|f\| var_*(\rho,\Delta).$$
\end{lem}
\begin{proof}
Let $(y,y')\in\delta$. By
construction, we have a bijection  $x\in
\pi^{-1}(y)\mapsto x'\in \pi^{-1}(y')$ such that
$(x,x')\in\Delta$. Then,
$$\begin{array}{lll}
E(f)(y)-E(f)(y')&={\displaystyle\sum_{\pi^{-1}(y)}} \rho
(x)[f(x)-f(x')]+\,[ \rho (x)- \rho(x')]f(x')\\
|E(f)(y)-E(f)(y')|&\le{\displaystyle\max_{\pi(x)=y}}|f(x)-f(x')|+\,\|f\|
{\displaystyle\max_{\pi(x)=y}}|{ \rho(x)\over
 \rho(x')}-1|\\
&\le var_+(f,\Delta)+\|f\| var_*(\rho,\Delta).
\end{array}$$
\end{proof}

Let us consider now  an AP equivalence relation $R$ on a compact
space
$X$ defined by a sequence
$$X_0\xrightarrow{\pi_{1,0}} X_1\rightarrow\ldots\rightarrow
X_{n-1}\xrightarrow{\pi_{n,n-1}}X_n
\rightarrow\ldots$$
where $X_0=X$ and for each $n\ge 1$, $X_n$ is a
Hausdorff space and $\pi_{n,n-1}$ is a surjective local
homeomorphism. We can construct inductively $\pi_{n,n-1}$-covers
$({\cal V}^n, {\cal U}^n)$ such
that ${\cal U}^n\prec {\cal V}^{n-1}$.

\begin{defn}\label{tower} Such a sequence $({\cal V}^n, {\cal U}^n)$ of 
$\pi_{n,n-1}$-covers will be
called a {\it tower} relative to the sequence $(X_n)$.
\end{defn}

Given $m<n$ and a sequence $a=(U_{m+1},U_{m+2},\ldots,U_n)$, where $U_k\in {\cal U}^k$,
such that $\,\overline U_{k+1}\subset
\pi_{k,k-1}(U_k)\,$ for $k=m+1,\ldots,n-1$, we define the following subset of $X_m$:
$$U_a=U_{m+1}\cap
\pi_{m+1,m}^{-1}(U_{m+2})\cap
\ldots\cap\pi_{n-1,m}^{-1}(U_n).$$
To avoid redundancies, we implicitly choose once for all, for $n=1,2,\ldots$ and
$U\in{\cal U}_n$ a set $V=V(U)\in{\cal V}_{n-1}$ such that $\overline U\subset V$.
We shall only consider sequences $\,a=(U_{m+1},U_{m+2},\ldots,U_n)\,$ with $U_k\in {\cal
U}^k$ such that $\pi_{k,k-1}(U_k)=V(U_{k+1})$ for
$m<k< n$.

\begin{lem}\label{tower} Let $({\cal V}^n, {\cal U}^n)$ be a tower. Given $m<n$, we
define
${\cal U}^{m,n}$ as the family of open sets
$U_a$, where 
$$a=(U_{m+1},\ldots,U_n),\quad U_k\in {\cal
U}^k,\quad \pi_{k,k-1}(U_k)=V(U_{k+1}),\quad  m<k< n.$$
Then $({\cal V}^n,{\cal
U}^{m,n}) $ is a $\pi_{n,m}$-cover.
\end{lem}
\begin{proof}
Let
$V\in {\cal V}^n$. Let us show that
$\pi_{n,m}^{-1}(V)$ is the disjoint union of the family of 
$U_a$'s, where $a=(U_{m+1},U_{m+2},\ldots,U_n)$ where $U_k\in {\cal U}^k$,
$\pi_{k,k-1}(U_k)=V(U_{k+1})$ for
$m<k< n$ and $V=\pi_{n,n-1}(U_n)$ and that
$\pi_{n,m}$ maps each $U_a$ homeomorphically onto
$V$. Let
$x_m\in
\pi_{n,m}^{-1}(V)$. We define
$x_k=\pi_{k,m}(x_m)$ for
$m<k\le n$. Since $x_n\in V$, there is a unique $\,U_n\in{\cal U}_n\,$ mapping
homeomorphically onto $V$ and containing  $x_{n-1}$. We define
$\,V_{n-1}=V(U_n)\in{\cal V}_{n-1}\,$ and
proceed by induction to construct $U_{n-1}\in{\cal U}_{n-1},\ldots, U_{m+1}\in{\cal
U}_{m+1}$. Then $x_m\in U_a$ where $\,a=(U_{m+1},U_{m+2},\ldots,U_n)$. If
$a=(U_{m+1},U_{m+2},\ldots,U_n)$ and
$a'=(U'_{m+1},U'_{m+2},\ldots,U'_n)$ are distinct, there exists a
larger
$k\le n$ such that $U_k\not=U'_k$. Then, $U_k$ and $U'_k$ are disjoint and so are
$U_a$ and
$U_{a'}$. The restriction of $\pi_{n,m}$ to $U_a$ is a composition of
homeorphisms and therefore maps $U_a$ homeomorphically onto
$V$.
\end{proof}
Thus, if we contract the initial sequence of spaces by means of a
subsequence $(n_k)$, our initial tower 
$({\cal V}^n,{\cal U}^n)$
provides a tower
$({\underline{\cal V}}^k,{\underline{\cal U}}^k)$ 
relative to
the sequence $({\underline X}_k=X_{n_k})$: we set ${\underline{\cal V}}^k={\cal
V}^{n_k}$ and ${\underline{\cal U}}^k={\cal U}^{n_{k-1},n_k}$.

In particular, we will denote ${\cal U}_n={\cal U}^{0,n}$.  Note that $({\cal
V}^n,{\cal U}_n)$ is a $\pi_n$-cover and that $\,{\cal U}_n\prec{\cal U}_{n-1}$.

Given a tower $({\cal V}^n, {\cal U}^n)$ as above for the defining
sequence $(X_n)$, we define the following entourages in $X\times X$:
$$\Delta_n=\Delta_{{\cal U}_n}=\cup\, \{ U\times U,\quad U\in {\cal U}_n\}$$
and in $X_n\times X_n$:
$$\delta_n=\Delta_{{\cal V}^n}=\cup\,\{ V\times V,\quad V\in {\cal V}^n\}.$$
Note that
$(\Delta_n)$ is a sequence of neighborhoods of the diagonal $\Delta_X$ of $X\times X$
such that
$\Delta_n\subset\overline\Delta_n\subset\Delta_{n-1}$.

\begin{defn}\label{generator} We say that the tower 
$({\cal V}^n, {\cal U}^n)$ for the defining
sequence $(X_n)$ is a {\it generator} if $\cap\Delta_n$ is reduced to the
diagonal $\Delta_X$ of $X\times X$. In other words, a tower $({\cal V}^n, {\cal U}^n)$
is a generator iff $(\Delta_n)$ is a fundamental system of neighborhoods of $\Delta_X$.
\end{defn}

\begin{ex}\label{Bratteli diagram}{\bf AF equivalence relation defined by a
Bratteli diagram}. Let $(V,E)$ be a Bratteli diagram, let  $X$ be its infinite
path space and
$X_n$ the space of paths starting at level $n$. The Bratteli diagram
defines a tower $({\cal V}^n, {\cal U}^n)$, where $\,{\cal V}^n\,$ is the
partition of $X_n$ by the cylinder sets 
$$Z^n(v)=\{x_{n+1}x_{n+2}\ldots\in X_n: s(x_{n+1})=v\},$$ where $v\in V(n)$ and
 $\,{\cal U}^n\,$ is the
partition of $X_{n-1}$ by the cylinder sets 
$$Z^{n-1}(e)=\{x_nx_{n+1}\ldots\in X_{n-1}: x_n=e\},$$ where $e\in
E(n)$. Then
$\,{\cal U}_n\,$ is the partition of $X$ by the cylinder sets $Z(a)$, where
$a=a_1\ldots a_n$ is a path from level $0$ to level $n$. We have, for
$x,y\in X$, 
$$(x,y)\in \Delta_n\Leftrightarrow x_1=y_1,\ldots, x_n=y_n$$
 and for
$x,y\in X_n$, 
$$(x,y)\in \delta_n\Leftrightarrow s(x)=s(y), $$
i.e $x$ and $y$ start
from the same vertex. 

Let $D\in
Z^1(R,{\bf R}_+^*)$ and let $(\rho_n)$ be the sequence of its normalized
potentials. We have
$$var_*(\rho_n,\Delta_n)=\sup |{\rho_n(ax)\over\rho_n(ay)}-1|$$
where the supremum is taken over all paths $a$ from level $0$ to level
$n$ and all infinite paths $x,y$ starting at $r(a)$. Note that if $D$ is a
quasi-product cocycle relative to $(V,E)$, then for all $n$ and for all paths $a$ from
level $0$ to level $n$, we have
$\displaystyle{\rho_n(ax)\over\rho_n(ay)}=1$ and therefore $var_*(\rho_n,\Delta_n)=0$.
The estimates of this section can be used to extend the results
of \corref{AF unique state} concerning quasi-product cocycles to cocycles for which
there is  a good  control of
$var_*(\rho_n,\Delta_n)$.
\end{ex}

Let us return to the general case of an AP equivalence relation $R$ with a given tower
$({\cal V}^n, {\cal U}^n)$ and let us consider the cocycles $D$ on $R$ which satisfy
$\lim var_*(\rho_n,\Delta_n)=0$, where the $\Delta_n=\Delta_{{\cal U}^n}$ and where
$(\rho_n)$ is a sequence of potentials for
$D$. It results from \lemref{normalization} that if this condition is satisfied by a
sequence of (unnormalized) potentials $(g_n)$, it is also satisfied by the sequence of
normalized potentials $(\rho_n)$. Let us also observe that, when the tower
$({\cal V}^n, {\cal U}^n)$ is a generator, this
condition is invariant under cohomology. Indeed, if
$$\underline D(x,y)=b(x)D(x,y)b(y)^{-1}$$
for some $b\in C(X,{\bf R}_+^*)$, it admits the sequence of potentials $(b\rho_n)$.
The condition $\lim var_*(\rho_n,\Delta_n)=0$ and the uniform continuity of $b$ imply
that $\lim var_*(b\rho_n,\Delta_n)=0$. Thus a cocycle
$D$ on an AF equivalence relation given by a Bratteli diagram $(V,E)$ which is
cohomologous to a quasi-product cocycle with respect to the diagram must satisfy
$\lim var_*(\rho_n,\Delta_n)=0$. One can prove a weak converse: if $\liminf
var_*(\rho_n,\Delta_n)=0$, then $D$ is cohomologous to a
quasi-product cocycle relative to a contraction of $(V,E)$. The proof is similar to that
of Theorem 3.1 of
\cite{ren:AF}.

\section{Stationary systems}

We shall now apply the above estimates to stationary systems,
in the following sense. We say that the sequence
$$X_0\xrightarrow{\pi_{1,0}} X_1\rightarrow\ldots\rightarrow
X_{n-1}\xrightarrow{\pi_{n,n-1}}X_n
\rightarrow\ldots$$
 is {\it stationary} if for all $n\in{\bf N}$, $X_n=X$ and for all $n\ge
1$,
$\pi_{n,n-1}=T$. As usual, we assume that $X$ is compact and that $T$ is
a local homeomorphism and a surjection. The AP equivalence
relation defined by this sequence is
 $$R=R(X,T)=\{(x,y)\in X\times X:\exists n\in{\bf N}: T^nx=T^ny\}.$$
Then
$R$ is a subgroupoid of the semi-direct product groupoid (see for example
\cite{ren:AF}) of the dynamical system 
$$G(X,T)=\{(x,m-n,y)\in X\times{\bf Z}\times X:m,n\in{\bf N},
T^mx=T^ny\}$$
We say that a cocycle $D\in Z^1(R,{\bf R}_+^*)$ is {\it stationary} if it
is the restriction of a cocycle in $Z^1(G(X,T),{\bf R}_+^*)$. Because
$G(X,T)$ is singly generated (see 4.1 of
\cite{ren:AF}),  cocycles in $Z^1(G(X,T),{\bf R}_+^*)$ are
in a one-to-one correspondence with functions
$g\in C(X,{\bf R}_+^*)$. However, the cocycles defined by $g$ and
$\lambda g$, where $\lambda\in{\bf R}_+^*)$, (or more generally
$\lambda\in C(X,{\bf R}_+^*)^R$) will have the same restriction to $R$.
More explicitly, a cocycle
$D\in Z^1(R,{\bf R}_+^*)$ is stationary iff it admits a sequence of 
potentials $(g_n)$ of the form:
$$g_n(x)=g(x)g(Tx)\ldots g(T^{n-1}x)$$
where $g$ is a given function in $C(X,{\bf R}_+^*)$.
Then the normalized potentials are given by 
$$\rho_n={g_n\over Z_n\circ T^n} ,\quad \hbox{where}\quad Z_n(x)=\sum_{T^ny=x}g_n(y)$$
are the partition functions.
The normalized local potentials are given by
$$\rho_{n,n-1}(x)=(Z_n(Tx))^{-1}g(x)Z_{n-1}(x).$$

We say that a tower $({\cal V}^n, {\cal U}^n)$  for the stationary
sequence 
\newline $(X_n=X,\pi_{n,n-1}=T)$ is stationary if it does not depend on
$n$. Thus, a stationary tower is given by a $T$-cover $({\cal V}, {\cal
U})$ such that ${\cal U}\prec{\cal V}$. As before, we
define ${\cal U}_n$ as the cover of $X$ by the open sets
$$U_a=U_1\cap
T^{-1}(U_2)\cap
\ldots\cap T^{-(n-1)}(U_n),$$
where $a=(U_1,\ldots, U_n)$, $U_1,\ldots, U_n\in {\cal U}$, $\overline
U_{k+1}\subset T(U_k)$. We also define the entourages
$\Delta_n=\Delta_{{\cal U}_n}\subset X\times X$ and
$\delta_n=\Delta_{{\cal V}^n}\subset X_n\times X_n$. Note that $\delta_n$
does not depend on $n$; we call it $\delta$. Remember that the tower is called a
generator if $\cap\Delta_n$ is reduced to the
diagonal $\Delta_X$ of $X\times X$. We
then say that the corresponding $T$-cover $({\cal V}, {\cal
U})$ is a generator. In the classical terminology of \cite{wal:introduction}, the
cover ${\cal U}$ itself is called a generator if $\cap\Delta_n$ is reduced to the
diagonal $\Delta_X$.

\begin{defn} Let $T:X\rightarrow X$ be a continuous map on a compact space.
One says that
$g\in C(X,{\bf R}_+^*)$ satisfies {\it Walters' condition} if for all $\epsilon>0$,
there exists an entourage $\Delta$ of the uniform structure of $X$ such that the
sequence
$$\Delta_n=\{(x,y)\in X\times X:\forall k=1,\ldots,n-1, (T^kx,T^ky)\in
\Delta\}$$ decreases to the diagonal and for all
$n$,
$var_*(g_n,\Delta_n)\le\epsilon$.
\end{defn}

An equivalent definition is that,
given a generator $({\cal V},{\cal U})$, for all $\epsilon$, there exists an integer $N$
such that for all $n$,
$$var_*(g_n,\Delta_{N+n})\le\epsilon.$$

We can make a similar definition for an arbitrary cocycle on an ``almost
stationary'' AP equivalence relation
$(X,R)$.

\begin{defn} Let $R$ be an AP equivalence relation on a compact space $X$ with a
defining sequence $(X_n,\pi_n)$ where $X_n=X$ for all $n$. We shall say that
$D\in Z^1(R,{\bf R}_+^*)$ satisfies {\it Walters' condition} if for all
$\epsilon>0$, there exists a generator $({\cal V}_n,{\cal U}_n)$ with ${\cal
V}_n={\cal V}_1$ for all $n\ge 1$ and a sequence of potentials $(g_n)$ for $D$ such
that for all
$n$,
$var_*(g_n,\Delta_{{\cal U}_n})\le\epsilon$
\end{defn}

\begin{thm}\label{Walters} Let $R$ be an AP equivalence relation on a compact
space $X$ with a defining sequence $(X_n,\pi_n)$ where $X_n=X$ for all $n$ and let
$D\in Z^1(R,{\bf R}_+^*)$. Assume that:
\begin{enumerate}
\item the AP equivalence relation $R$ is minimal;
\item $D$ satisfies Walters'condition.
\end{enumerate}
Then, $D$ is uniquely ergodic. 
\end{thm}

\begin{proof} The proof is essentially the same as in Theorem 6 of \cite{wal:78}
(see also \cite{zin:thermo}). We introduce the normalized potentials
$(\rho_n)$ and their expectations
$(E_n)$. According to \lemref{normalization}, they also satisfy  Walters'
condition. Let
$\epsilon>0$ be given; there exists a generator $({\cal V}_n, {\cal U}_n)$ such
that for all
$n\ge 1$, ${\cal V}_n={\cal V}_1$ and
$var_*(\rho_n,\Delta_n)\le\epsilon$. We apply \lemref{Walters} with $\pi=\pi_n$ and
$({\cal V}_1,{\cal U}_n)$: given $f\in C(X)$, we have
$$var_+(E_n(f),\delta)\le var_+(f,\Delta_n)+\|f\| var_*(\rho_n,\Delta_n),$$
where $\delta=\Delta_{{\cal V}_1}$ is an entourage of the uniform structure of $X$.
This show the equicontinuity of the sequence $(E_n(f))$ in
$C(X)$. Since this sequence is also bounded, there exists a sequence $(n_k)$
tending to infinity and
$f^*\in C(X)$ such that $E_{n_k}(f)$ converges uniformly to $f^*$. Then, for any $m$,
$E_{m+n_k}(f)$ converges uniformly to $E_m(f^*)$.
Since the sequence $(E_n(f)_{max})$ is decreasing, it converges to
$f^*_{max}=E_m(f^*)_{max}$. Let $S_m$ be the set of points
where $E_m(f^*)$ takes its maximum. Note that, because of the strict positivity of
$\rho_{m+k,m}$, $\pi_{m+k,m}^{-1}(S_{m+k})\subset S_m$. Therefore
$(\pi_m^{-1}(S_m))$ is a decreasing sequence of non-empty closed sets and has a
non-empty closed intersection $S$. Moreover $S$ is invariant under $R$. Because
of $(i)$, $S=X$ and $f^*$ is a constant function. This implies that
$var_+(E_{n_k}(f))$ tends to zero. Since the sequence $(var(E_n(f))$ is decreasing, it
also tends to zero. Then \corref{unique ergodicity} gives the unique ergodicity.
\end{proof}

\begin{rem} Under the assumptions of the theorem,
let $\mu$ be the unique quasi-invariant measure on $X$ admitting $D$ as
Radon-Nikodym derivative. Then, for $f\in C(X)$, $(E_n(f))$ converges
uniformly to $\mu(f)1_X$. In particular, $(E_n(f)(x))$ converges
to $\mu(f)$ for all $x\in X$.

\end{rem}

In the case of a stationary cocycle, one retrieves the well-known
uniqueness result:

\begin{cor}\label{Theorem 6} (Theorem 6,\cite{wal:78}) Let $X$ be a compact space,
let
$T:X\rightarrow X$ be a surjective local homeomorphism and let
$g\in C(X,{\bf R}_+^*)$. Assume that:
\begin{enumerate}
\item the AP equivalence relation $R(X,T)$ is minimal;
\item there
exists an integer $L\ge 1$ and a
$T^L$-cover
$({\cal V}, {\cal U})$ which is a generator in the above sense for the
sequence defined by
$T^L$;
\item there
exists an integer $M\ge 1$ such that $g_M$ satisfies Walters' condition with
respect to $T^M$.
\end{enumerate}
Then, the stationary cocycle $D$ defined by $g$
on $R(X,T)$ is uniquely ergodic. 
\end{cor}

\begin{proof} Let $N=LM$. If
$({\cal V},{\cal U})$ is a generator for $T^L$, then
$({\cal V},{\cal U}_M)$ is a generator for $T^N$. If $g_M$ satisfies Walters'
condition with respect to $T^M$, $g_N$ satisfies  Walters' condition with
respect to $T^N$. Moreover
$R(X,T^N)=R(X,T)$. Therefore, replacing $T$ by
$T^N$, we may assume that $L=M=1$. Then, the AP equivalence relation
$R(X,T)$ defined by $(X,T)$ and the stationary cocycle $D\in Z^1(R,{\bf
R}_+^*)$ satisfy the assumptions of the theorem.
\end{proof}

\begin{ex}\label{Subshifts of finite type} {\it Subshifts of finite type}. Let
$\Gamma=(\Gamma^{(0)},\Gamma^{(1)})$ be a finite graph. We assume that each vertex
receives and emits at least one edge. We let
$X$ be the space of one-sided infinite paths $x=x_1x_2\ldots$ and
$T:X\rightarrow X$ be the one-sided shift $T(x_1x_2\ldots)=x_2\ldots$. For
$v\in\Gamma^{(0)}$, we let $Z(v)$ be the set of paths starting at $v$ and for
$e\in\Gamma^{(1)}$, we let $Z(e)$ be the set of paths having $e$ as initial edge.
We define ${\cal V}=\{Z(v),  v\in\Gamma^{(0)}\}$ and ${\cal U}=\{Z(e), e\in
\Gamma^{(1)}\}$. Then $({\cal V}, {\cal U})$ is a generator for the
sequence defined by $T$. This is a particular case of \exref{Bratteli diagram} where
$(V,E)$ is the stationary Bratteli diagram defined by $\Gamma$. It is known that
$R(X,T)$ is minimal if and only if the graph
$\Gamma$ is primitive (i.e. there exists an integer $L\ge 1$ such that every pair of
vertices $(v,w)$ can be joined by a path of length $L$. This is also
equivalent to $T$ being topologically mixing. Walters' condition for
$g\in C(X,{\bf R}_+^*)$ reads here:
$$\forall \epsilon>0,\, \exists N\in{\bf N}:\,   x_1\ldots
x_{n+N}=y_1\ldots y_{n+N}\Rightarrow |{g_n(x)\over g_n(y)}-1|\le\epsilon.$$
This is also a particular case of the next example.
\end{ex}
\begin{ex}\label{Expansive maps} {\it Expansive maps}. Let $X$ be a compact metric
space and let $T:X\rightarrow X$ be surjective and a local homeomorphism. One says
that $T$ is (positively) expansive if there exists $\epsilon>0$ such that for
every pair of points of $X$,
$x\not=y$, there is $n\in{\bf N}$ such that
$d(T^nx,T^ny)\ge\epsilon$.
\begin{lem} Let $T:X\rightarrow X$ be a surjective and expansive local homeomorphism.
Then, there exist an integer $L\in{\bf N}$ and a $T^L$-cover $({\cal V},{\cal
U})$ which is a generator for $T^L$.
\end{lem}
\begin{proof}
 Replacing the metric by a topologically equivalent metric, we may assume
that $T$ is locally expanding: there exists
$\tau>0$ and
$\lambda<1$ such that  $d(x,y)<\tau\,\Rightarrow\, d(x,y)\le
\lambda d(Tx,Ty)$. Let $({\cal V},{\cal U})$ be a finite $T$-cover such that the
elements of $\cal U$ have diameter strictly less than $\tau$. Let $c$ be an upper
bound for the diameter of the elements of $\cal V$. For each
$L\in{\bf N}$,
$({\cal V},{\cal U}_L)$ is a $T^L$-cover. For $L$ sufficiently large, the diameter of
each element of ${\cal U}_L)$, which is majorized by $\lambda^lc $, will be strictly
less than the Lebesgue number of the cover $\cal V$ and we will have ${\cal
U}_L\prec {\cal V}$. Let $\Delta_n=\Delta_{{\cal U}_n}$. Since
$d(\Delta_n)\le
\lambda^nc$ tends to $0$, $({\cal V},{\cal U}_L)$ is a generator.
\end{proof}

One says that $T:X\rightarrow X$ is {\it exact} if for any
non-empty open set
$U$, there is an integer $n>0$ such that $T^n(U)=X$. This condition is
equivalent to the minimality of
$R(X,T)$ (for equivalent conditions, see \cite{kr:kms}).

Let $g\in C(X,{\bf
R}_+^*$ and let $(g_n)$ be the corresponding sequence of potentials. We fix a
compatible metric
$d$ on
$X$ and define 
$$d_n(x,y)=\max\{d(Tx,Ty),\ldots,d(T^{n-1}x,T^{n-1}y)\}.$$ The usual
Walters condition for $g$ is that
for all $\epsilon>0$,
there exists $\delta>0$ such that for all
$n$,
$var_*(g_n,\Delta_n)\le\epsilon$, where 
$$\Delta_n=\{(x,y)\in X\times X: d_n(x,y)<\delta\}.$$
(This condition is often expressed in terms of $\varphi=\log g$.) Since
the sets
$$\Delta=\{(x,y)\in X\times X: d(x,y)<\delta\}$$
form a fundamental system of entourages of the uniform structure of $X$, this
condition is equivalent to ours. Thus one retrieves the well-known result that if 
$T$ is expansive and exact and if for some integer $M\ge 1$, $g_M$ satisfies
Walters' condition with respect to $T^M$,
then the stationary cocycle $D$ defined by $g$
on $R(X,T)$ is uniquely ergodic. 
\end{ex}

The above theorem establishes the unique ergodicity of quasi-product
cocycles on stationary Bratteli diagrams:

\begin{cor} Let $(V,E)$ be a stationary Bratteli diagram defined by a a
finite primitive graph. Let $X=X(V,E)$ and $R=R(V,E)$. Then every
quasi-product cocycle
$D\in Z^1(X,R)$ is uniquely ergodic.
\end{cor}

\begin{proof}
As we have seen, $R(V,E)$ is minimal. The condition $\lim
var_*(\rho_n,\Delta_n)=0$, where we use the notation of \exref{Bratteli
diagram}, is trivially satisfied for a quasi-product cocycle and implies
Walters'condition for
$D$.
\end{proof}

\begin{rem}
This result shows that the sufficient condition $(ii)$ of \corref{AF
unique state} is
not necessary.
\end{rem}

\section{Application to the transfer operators}

We show in this section how our above results on the unique ergodicity of
cocycles on AP equivalence relations imply parts of the
Perron-Frobenius-Ruelle theorem (as stated for example in
\cite{wal:78}) namely the existence and the uniqueness of the Perron eigenvalue
and eigenvector of the Ruelle operator
$$\,{\cal L}^*_g\,:C(X)^*\rightarrow C(X)^*.$$ 
The setting is a dynamical system $(X,T)$, where $X$ is a compact space and
$T$ is a  surjective local homeomorphism of $X$ onto itself.

An arbitrary continuous function $g\in C(X,{\bf R}_+^*)$ defines both a cocycle
$D_g$ in $Z^1(G,{\bf R}_+^*)$, where the definition of $G=G(X,T)$ has been recalled
earlier, and a {\it Ruelle} (or {\it transfer}) operator ${\cal L}_g:
C(X)\rightarrow C(X)$ according to
$${\cal L}_g(f)(x)=\sum_{Ty=x}g(y)f(y).$$
The relation between the Ruelle operator ${\cal L}_g$ and the cocycle $D_g$ is
given by the following elementary result.
\begin{prop} (Proposition 4.2 of \cite{ren:AF}) Let $\mu$ be a probability measure
on $X$. The following conditions are equivalent:
\begin{enumerate}
\item $\mu$ is quasi-invariant with respect to $G(X,T)$ and admits $D_g$ as
Radon-Nikodym derivative;
\item $\, {\cal L}^*_g\,\mu=\mu$.
\end{enumerate}
\end{prop}
The quasi-invariance with respect to $G(X,T)$ is defined in the same way that the
quasi-invariance with respect to $R(X,T)$: it means that the measures $r^*\mu$ and
$s^*\mu$ are equivalent, where $r,s:G(X,T)\rightarrow X$ are the projections.
Since $R(X,T)\subset G(X,T)$, the quasi-invariance with respect to $G(X,T)$ implies
the quasi-invariance with respect to $R(X,T)$ and the restriction of a
Radon-Nikodym derivative relative to $G(X,T)$ is a
Radon-Nikodym derivative relative to $R(X,T)$.

\begin{prop}\label{Ruelle} Let $(X,T)$ be a dynamical system as above and let $g\in
C(X,{\bf R}_+^*)$. Let $G(X,T), R(X,T), D_g$ and ${\cal L}_g$ be as above. Assume
that the restriction $D$ of $D_g$ to $R(X,T)$ is uniquely ergodic. Then the
eigenvalue problem 
$$ {\cal L}^*_g\,\mu=\lambda\mu , $$
where $\mu$ is a probability measure, admits one and only one
solution $\mu$. 
\end{prop}

\begin{proof}
Let us first show the uniqueness of $\mu$.
Suppose that the probability measure $\mu$ is an eigenvector of the transpose of
the Ruelle operator:
$ {\cal L}^*_g\,\mu=\lambda\mu$. Then $\lambda>0$ and 
$ ^t{\mathcal L}_{\lambda^{-1}g}\mu=\mu$. Therefore, $\mu$ is quasi-invariant with
respect to
$G(X,T)$ with Radon-Nikodym derivative
$D_{\lambda^{-1}g}$. Since
$D_{\lambda^{-1}g}$ and $D_g$ have the same restriction $D$ to $R(X,T)$, $\mu$ is
solution of the Radon-Nikodym problem for $D$. This shows that $\mu$, if it
exists, is unique.

Let us show that the unique solution $\mu$ of the Radon-Nikodym problem for $D$ is
a solution of the eigenvalue problem. The cocycle $D$ on the AP equivalence
relation $R(X,T)$ admits the sequence $(g_n)$ as (unnormalized) potentials. The
corresponding local potentials are $g_{n,n-1}=g$ for all $n$. In other words,
according to \propref{state space} and the Appendix, the
solutions of the Radon-Nikodym problem for $D$ are in one-to-one
correspondance with the states of the inductive limit of the sequence
$$C(X)\xrightarrow{{\mathcal L}_g} C(X)\rightarrow\ldots\rightarrow
C(X)\xrightarrow{{\mathcal L}_g}C(X)
\rightarrow\ldots$$
Explicitly, a state is given by a sequence of measures $(\mu_n)$ such that $\mu_0$
is a probability measure and $\mu_n= \,{\cal L}^*_g\,\mu_{n+1}$ for all $n\in{\bf
N}$. The correspondence is $(\mu_n)\mapsto \mu_0$. It is then clear that, if
$(\mu_n)$ is a state, so is $(\mu'_n=\,{\cal L}^*_g\,\mu_n/\lambda)$, where
$\lambda=\mu_0({\mathcal L}_g(1))$. By unique ergodicity of $D$, we obtain that
the above measure $\mu$ satisfies
$\mu=\,{\cal L}^*_g\,\mu/\lambda$, where $\lambda=\mu_0({\mathcal L}_g(1))$.

\end{proof}

\appendix

\section{State space of dimension groups}

Since the work \cite{ell:ind} of G.~Elliott, dimension groups have been used as a
convenient tool in the study of AF-algebras and topological Markov chains
(\cite{kri:dimension}). Our focus will be the state space of the dimension group,
which describes the traces of the AF-algebra. More generally, we shall use
dimension groups in the study of Radon-Nikodym cocycles on AP relations. Most
results of this section are not new. Those concerning dimension groups can be found
in the monograph \cite{goo:poag} by K.~Goodearl, in the general setting of
(partially)ordered abelian groups. Another basic reference on dimension groups is
\cite{eff:dimensions} by E.~Effros. We also give a slight improvement of a result
of A.~T\"or\"ok \cite{toe:AF} on the uniqueness of a trace on an AF-algebra.
The purpose of this appendix is to present the results we need, in the setting the
most appropriate to the Radon-Nikodym problem. We first introduce some notation.

When $X$ is a compact space, $C(X)$ designates
the real vector space of real-valued continuous functions on $X$ endowed
with the uniform norm $\|f\|=\sup_X|f(x)|$. We write $f\ge 0$ if $f(x)\ge 0$ for
all $x\in X$ and we denote by $1=1_X$ the constant function $1(x)=1$. 
This turns $(C(X),1_X)$ into an ordered real vector space with order-unit. We write
$f>0$ if $f(x)>0$ for all $x\in X$.

\begin{defn}\label{positive} A bounded linear operator
$A:C(X)\rightarrow C(Y)$, where $X,Y$ are compact spaces is called
\begin{itemize}
\item {\it positive} if $f\ge 0\Rightarrow Af\ge 0$;
\item {\it strongly positive} if $f > 0\Rightarrow Af > 0$;
\item {\it Markovian} if it is positive and $A1_X=1_Y$.
\end{itemize}
\end{defn}

When $X=\{1,\ldots,q\}$ and $Y=\{1,\ldots,p\}$,
$C(X)={\bf R}^q, C(Y)={\bf R}^p$ and $A$ is given by a matrix
$(a(i,j))\in M_{p,q}({\bf R})$. Positivity means $a(i,j)\ge 0$ for all
$i,j$ and strong positivity means positivity and non-zero rows; it is a weaker
condition than the strict positivity of the $a(i,j)$'s.

Let $(X_n), n\in{\bf N}$ be
a sequence of compact spaces
and let $\underline A=(A_n:C(X_{n-1})\rightarrow C(X_n)), n=1,2,\ldots$ be a
sequence of strongly positive operators. Our problem is to study the set of
sequences $\mu=(\mu_n)$, where
$\mu_n$ is a positive linear functional on $C(X_n)$, i.e. a measure on
$X_n$, satisfying the recurrence relation $\mu_{n-1}=A_n^*\mu_n$ for all
integers $n\ge 1$ as well as the normalization condition
$\mu_0(1_{X_0})=1$.

We introduce
$$u_n=A_n A_{n-1}\ldots A_11_{X_0}.$$
We view $({\cal E}_n=C(X_n),u_n)$ as an ordered real vector space with order-unit.
The sequence $\underline A$ defines an inductive system  of ordered real vector
spaces with order-unit and we can consider its inductive limit ${\cal E}={\cal E}(\underline
A)=\indlim {\cal E}_n$. It is an ordered real vector spaces with order-unit $u$. We denote
by
$j_n:{\cal E}_n\rightarrow {\cal E}$ the canonical morphisms. 

Transposition gives the
projective system
$(A^*_n:C(X_{n-1})^*\leftarrow C(X_n)^*), n=1,2,\ldots$.

We introduce the state space of $({\cal E}_n,u_n)$:
$$S_n=\{\tau\in C(X_n)^*: \hbox { positive and }\,\tau(u_n)=1\}.$$
It is a convex compact subset of $C(X_n)^*$ in the $*$-weak topology, in
fact it is a Choquet simplex. We are interested in its projective limit:
$$S=S(\underline A)=\projlim S_n=\{(\tau_n): \tau_n\in
S_n;\tau_n=\tau_{n+1}A_{n+1},\quad n\in{\bf N}\}.$$

It is known
(\cite{goo:poag}, 6.14) that $S$ is the state space of $({\cal E},u)$, i.e. the convex
set of positive homomorphisms $\tau: {\cal E}\rightarrow{\bf R}$ such that
$\tau(u)=1$. The value of
$\tau=(\tau_n)\in S$ on $f=j_n(f_n)\in {\cal E}$ is $\tau(f)=\tau_n(f_n)$. It is also
known (\cite{goo:poag}, 10.21) that 
$S$ is a (non-empty!) convex compact set,
and a Choquet simplex .

\begin{defn}\label{state space} We call ${\cal E}={\cal E}(\underline
A)=\indlim {\cal E}_n$ the dimension group and $S=S(\underline A)$ the state space of the
sequence
$\underline A=(A_n)$.
\end{defn}

\begin{rem} When $\underline A=(A_n)$ is a sequence of matrices with coefficients
in ${\bf N}$, the dimension group
${\cal E}(\underline A)=\indlim(A_n:C(X_{n-1},{\bf R})\rightarrow C(X_n,{\bf R}))$
is a coarser invariant than the usual dimension group
$K_0(\underline A)=\indlim(A_n:C(X_{n-1},{\bf Z})\rightarrow C(X_n,{\bf Z}))$
introduced by Elliott. More precisely ${\cal E}(\underline A)$ is the realification  of
$K_0(\underline A)$. However, the state spaces are the same.
\end{rem}

Let $\underline A =(A_n:C(X_{n-1})\rightarrow C(X_n))$ and
$\underline B=(B_n:C(Y_{n-1})\rightarrow C(Y_n))$ be two sequences of strongly
positive operators as above. A morphism from $\underline B$ to $\underline A$ is a
sequence
$\underline D=(D_n)$ of positive operators
$D_n:C(Y_n)\rightarrow C(X_n)$ such that $D_0(1_{Y_0})=1_{X_0}$ and which
makes the following diagram commutative:

\begin{equation*}
\begin{CD}
C(X_0)   @>A_1>>   \ldots C(X_{n-1})@>A_n>>C(X_n)\ldots\\
@AA{D_0}A          @AA{D_{n-1}}A             @AA{D_n}A\\
C(Y_0)   @>B_1>>   \ldots C(Y_{n-1})@>B_n>>C(Y_n)\ldots
\end{CD}
\end{equation*}

This implies that the $D_n$'s are strongly positive.
It induces a morphism $D:{\cal E}(\underline B)\rightarrow {\cal E}(\underline A)$ and a
morphism of their state spaces
$D^*:S(\underline A)\rightarrow S(\underline B)$ which sends $(\tau_n)$ into
$(\tau_nD_n)$. 

\begin{defn}\label{contraction} A {\it contraction} of a sequence of
strongly positive operators  $A_n:C(X_{n-1})\rightarrow C(X_n)$ is a sequence
$(B_k=A_{n_k}A_{n_k-1}\ldots A_{n_{(k-1)}+1})$, where $(n_k)$ is a strictly
increasing sequence of integers.
\end{defn} 

Note that $(A_n)$ and a contraction $(B_k)$ give the same inductive limit
${\cal E}$ and the same state space. The following elementary observation
reduces the problem to the Markovian case.

\begin{prop}\label{Markov reduction} Let $(A_n)$ be a sequence of strongly positive
operators. Then, there exists an isomorphic sequence consisting of
Markovian operators.
\end{prop}

\begin{proof} Let $D_n$ be the operator of multiplication by
$u_n$ and define
$$B_n=D_n^{-1}A_nD_{n-1}.$$
\end{proof}

In the study of the inductive limit and its state space, we may
therefore assume that we have a sequence $\underline B=(B_n)$ of Markovian
operators. That is what we do from now on.

The following result is well known in the theory of dimension
groups (see \cite{eff:dimensions}, Corollary 4.2. and \cite{goo:poag}, Corollary
4.10), where its proof is based on an ordered group analogue of the Hahn-Banach
theorem due to Goodearl and Handelman). We use here a compactness argument
modelled after
\cite{wal:Bowen}, Theorem 2.9.  Albeit elementary, it is one of the most useful
tools to compute the state space of a dimension group.

\begin{lem}\label{key lemma} Let $\underline B=(B_n)$ be a sequence of Markovian
operators, ${\cal E}={\cal E}(\underline B)$ and $S=S(\underline B)$ as above. Let
$f=j_m(f_m)\in {\cal E}$ and let $f_n=B_nB_{n-1}\ldots B_{m+1}f_m$ for $n\ge
m+1$. The following conditions are equivalent:
\begin{enumerate}
\item the sequence  $(\|f_n\|)$ tends to $0$,
\item $f$ belongs to $S_\perp$, i.e. $\mu(f)=0$ for all $\mu\in S$.
\end{enumerate}
\end{lem}

\begin{proof} $(i)\Rightarrow(ii)$ Given $\epsilon>0$, there exists $n$
such that $\|f_n\|<\epsilon$. Therefore,
$$|\mu(f)|=|\mu_n(f_n)|\le\|f_n\|<\epsilon$$ and
$\mu(f)=0$.

$(ii)\Rightarrow(i)$ Suppose that $(i)$ is not satisfied. Since the
sequence $(\|f_n\|), n=m, m+1,\ldots,$ is
decreasing, there exists $\epsilon>0$ such that
$\|f_n\|\ge\epsilon$ for all $n\ge m$. Therefore, for $n\ge m$, the
set $K_n=\{\mu_n\in S_n: |\mu_n(f_n)|\ge\epsilon\}$ is not empty.
Consider the sequence of subsets of $\prod_0^\infty S_n$:
$$S^N=\{(\mu_n):\mu_{n-1}=\mu_nB_n\quad{\rm for}\quad
n\le N\quad{\rm and}\quad \mu_n\in K_n\quad{\rm for}\quad
n\ge N\}$$
where $N\ge m$. It is a decreasing sequence of non-empty compact sets.
Then, an element $\mu=(\mu_n)$ of its intersection belongs to $S$ and
$\mu(f)\not=0$ since $|\mu(f)|=|\mu_m(f_m)|\ge\epsilon$.
\end{proof}

Let us first apply this lemma to study when the inductive system $\underline
A=(A_n)$ has a {\it unique state}, that is, when $S(\underline A)$ is
reduced to one element.

Given $f\in C(X)$, where $X$ is a compact space, we define:
$$f_{min}=\min_Xf(x),\quad  f_{max}=\max_Xf(x),\quad {\rm
var}(f)=f_{max}- f_{min}.$$

\begin{cor}\label{CNS unique state} Let $\underline B=(B_n)$ be a sequence of
Markovian operators. The following conditions are equivalent:
\begin{enumerate}
\item For all $m\in{\bf N}$ and all $f\in C(X_m)$, ${\rm
var}(B_nB_{n-1}\ldots B_{m+2}B_{m+1}f)$ tends to zero as $n$ tends to
infinity.
\item For all $m\in{\bf N}$ and all $f\in C(X_m)$, there exists
$\mu_m(f)\in{\bf R}$ such that
$\|B_nB_{n-1}\ldots B_{m+2}B_{m+1}f-\mu_m(f)1_{X_n}\|$ tends to zero as $n$ tends to
infinity.
\item The state space $S(\underline B)$ defined above
consists of a single element. If $\nu=(\nu_n)$ is this element, then for
all
$m\in{\bf N}$,
$\nu_m$ agrees with $\mu_m$ defined in
$(ii)$.
\end{enumerate}
\end{cor}

\begin{proof} Given $f_m\in C(X_m)$, we define $f_n=B_nB_{n-1}\ldots
B_{m+1}f$ for $n\ge m+1$. Let us show that $(i)\Rightarrow (ii)$. More
precisely, we show that the condition $(i)$ with $m=0$ implies the
condition $(ii)$ with
$m=0$. Let
$f\in C(X_0)$. The sequences
$(f_n^{min})$ and $(f_n^{max})$ are adjacent. If moreover
$f_n^{max}-f_n^{min}$ tends to zero, they converge to the same limit
$\mu(f)\in{\bf R}$. The inequality
$$\|f_n-\mu(f)1_{X_n}\|\le{\rm var}(f_n)$$
gives the conclusion. The reverse implication $(ii)\Rightarrow (i)$ is
clear.

Let us show that $(ii)\Rightarrow (iii)$. We know that $S(\underline B)$ is not
empty. Let $(\nu_n)\in S(B)$. We fix $m\in{\bf N}$. For all $n\in{\bf
N}, n\ge m+1$, we can write
$\nu_m=\nu_nB_nB_{n-1}\ldots B_{m+1}$. Therefore for all $f\in C(X_m)$,
$$\nu_m(f)-\mu_m(f)=\nu_n(B_n\ldots B_{m+1}f-\mu_m(f)1_{X_n}),$$
$$|\mu_0(f)-\mu(f)|\le \|B_n\ldots B_{m+1}f-\mu_m(f)1_{X_n}\|.$$
One concludes that $\nu_m(f)=\mu_m(f)$ and $\nu_m=\mu_m$.

Let us show that $(iii)\Rightarrow (ii)$. Let $\nu=(\nu_n)$ be the
single element of $S$. For $f_m\in C(X_m)$, $j_m(f_m-\nu_m(f_m)1_{X_m})$
belongs to $S_\perp$. According to the lemma, this implies $(ii)$.
\end{proof}

Thus the above corollary gives a necessary and sufficient condition for an
inductive system $\underline A=(A_n)$ to
have a unique state (other necessary and sufficient conditions can be found
in \cite{goo:poag}, Chapter 4, in the general setting of partially ordered abelian
groups). We shall give a more practical sufficient condition, along the lines of
the work \cite{toe:AF} of A.~T\"or\"ok, who studied the uniqueness of traces on
AF C$^*$-algebras. We first replace the sequence
$\underline A=(A_n)$ by the Markovian sequence
$\underline B=(B_n)$, where $B_n=D_n^{-1}A_nD_{n-1}$ and $D_n$ is the operator of
multiplication by
$u_n=A_n A_{n-1}\ldots A_11_{X_0}$ in order to apply $(i)$ of \corref{CNS unique
state}. It will be used under the following form.

\begin{cor} Let $\underline B=(B_n)$ be a sequence of Markovian
operators. Suppose that there exists a sequence of positive numbers $(\epsilon_n)$
such that $\sum\epsilon_n=\infty$ and ${\rm var}(B_nf)\le (1-\epsilon_n){\rm
var}(f)$ for all $n$ and all $f\in C(X_{n-1})$. Then $\underline B$ has a unique
state.
\end{cor}
\begin{proof} By induction, for all $m$, all $f\in C(X_m)$ and all $n\ge m+1$,
$${\rm var}(B_nB_{n-1}\ldots B_{m+1}f)\le\prod_{m+1}^n (1-\epsilon_k){\rm
var}(f).$$
This tends to $0$ when $n$ tends to $\infty$. 
\end{proof}

It remains to establish an estimate on the variation ${\rm var}(B_nf)$. We shall
only consider here the case when the $X_n$'s are finite
sets (i.e. the $B_n$'s are Markovian matrices). Let us start with the following
lemma.

\begin{lem}\label{torok}(cf. \cite{toe:AF}, Lemma 2.) Let $I,J$ be finite sets and
let
$B:C(I)\rightarrow C(J)$ be a Markovian operator defined by a Markovian matrix
$b:J\times I\rightarrow{\bf R}$. Let
$$\epsilon=\min\{\sum_{i\in I_1}b(j,i)+\sum_{i\notin I_1}b(j',i): j,j'\in J,\quad
 I_1\subset I\}.$$
Then, for all $f\in C(I)$, ${\rm var}(Bf)\le (1-\epsilon){\rm var}(f)$ .
\end{lem}

\begin{proof} We first note that ${\rm var}:C(I)\rightarrow {\bf R}_+$ is a
semi-norm. It is invariant under translation by
${\bf R}1_I$ and induces a norm on the quotient space $C(I)/{\bf R}1_I$.
If
${\rm var}(f)=0$,
$f=c1_I$, where
$c\in{\bf R}$,
$Bf=c1_J$ and 
${\rm var}(Bf)=0$. The case ${\rm var}(f)\not=0$ can be reduced to the case
$f_{min}=0, f_{max}=1$ by considering $g=(f-f_{min}1_I)/{\rm var}(f)$. We will get our
estimate if we show that $var(Bf)\le (1-\epsilon)$ for all $f$ in the convex set
$$C=\{f\in C(I):\quad \forall i\in I, 0\le f(i)\le 1\}.$$
Since the map $f\mapsto var(Bf)$ is convex, its maximum is attained on the set
$\partial C$ of extremal points of $C$:
$$\max_C var(Bf)= \max_{\partial C} var(Bf)=\max_{I_1\subset I} var (B1_{I_1}).$$
For $I_1\subset I$ and $j,j'\in J$, we have
\begin{eqnarray}
B1_{I_1}(j)-B1_{I_1}(j')=&\sum_{i\in I_1}b(j,i)-\sum_{i\in I_1}b(j',i)\nonumber\\
=&1-\sum_{i\notin I_1}b(j,i)-\sum_{i\in I_1}b(j',i)\nonumber
\end{eqnarray}
Since $\sum_{i\notin I_1}b(j,i)+\sum_{i\in I_1}b(j',i)\ge\epsilon$, $var (B1_{I_1})\le
1-\epsilon$.
\end{proof}

\begin{rem} This lemma gives, in the case $I=J$, an estimate of the spectral gap
of the Markovian matrix $B$ and a Perron-Frobenius theorem for such a matrix. As
noted in \cite{toe:AF}, the matrix need not be primitive: a necessary and
sufficient condition for $\epsilon>0$ is that $B$ does not have any pair of
orthogonal rows.
\end{rem}

We deduce from the above lemma the following result, valid for strongly positive
matrices rather than Markovian matrices, which will be used to show the convergence
of stationary systems.

\begin{lem}\label{key estimate} Let $J,I$ be finite sets and let $A:C(I)\rightarrow
C(J)$ be an operator defined by a strongly positive matrix $a:J\times
I\rightarrow{\bf R}_+^*$. Let $u\in C(I)$ be a strictly positive vector
and let 
$B:C(I)\rightarrow C(J)$ be the
Markovian operator  defined by the matrix
 $$\,b(j,i)=(Au(j))^{-1}a(j,i)u(i).$$
Then, for all $f\in C(I)$, $${\rm var}(Bf)\le \big(1- {a_{min}\over
a_{max}}\big)\,{\rm var}(f).$$
\end{lem}

\begin{proof} It suffices to show that the constant $\epsilon$ of the previous
lemma satisfies $\epsilon\ge a_{min}/a_{max}$. This is immediate, because of
the obvious inequalities 
$$\sum_{i\in I_1}a(j,i)u(i)\ge a_{min}\sum_{i\in I_1}u(i)\quad {\rm and}\quad
Au(j)\le a_{max}\sum_{i\in I}u(i).$$
which are valid for all $I_1\subset I$ and $j\in J$.
\end{proof}

Thus one gets a convenient condition on the sequence $(A_n)$ ensuring unique state.
A similar result is given in  \cite{toe:AF} but ours has the advantage to be
expressed directly in terms of the $(A_n)$'s rather than the $(B_n)$'s.

\begin{cor}\label{CS unique state} Let $\underline A =(A_n:C(I_{n-1})\rightarrow
C(I_n))$  be a sequence of strongly positive operators, where the $I_n$'s are
finite spaces. Let
$\epsilon_n$ be the ratio of the smallest matrix coefficient of $A_n$ over its
largest matrix coefficient. If $\sum\epsilon_n=\infty$, then $\underline A$ has a
unique state.
\end{cor}

As an application, let us study the case of a stationary sequence $(A_n=A)$.

\begin{ex}\label{Perron-Frobenius} {\it The Perron-Frobenius theorem for primitive
matrices}.  Let $I$ be a
finite set and let $a:I\times I\rightarrow {\bf R}_+$ be a primitive matrix; this
means that there exists a positive integer $L$ such that $a^L$ is strictly
positive. The matrix $a$ defines a strongly positive linear operator
$A:C(I)\rightarrow C(I)$. Let us first show that the stationary sequence
$\underline A=(A_n=A)$ has a unique state. Since contracting does not change the
state space, we consider instead the stationary sequence
$\underline A'=(A'_n=A^L)$. Using above notation,
$\epsilon_n=\epsilon=a^L_{min}/a^L_{max}$ is constant and strictly positive;
therefore the condition $\sum\epsilon_n=\infty$ is satisfied.  This shows that
$S(\underline A)$ has a unique element $\mu=(\mu_n)$. Recall that this is a
sequence of (positive) measures on $I$ satisfying $\mu_{n-1}=A^*\mu_n$ and
$\mu_0(1)=1$. Since the sequence $\nu=(\nu_n)$ where $\nu_n=\lambda\mu_{n+1}$ and
$\lambda=1/\mu_1(1)$ is also a state, we must have $\mu=\nu$. This shows that
$\mu_0$ is an eigenvector of
$A^*$ for the eigenvalue $\lambda$ and that $\mu_n=\lambda^{-n}\mu_0$. Conversely,
if $\nu_0$ is a probability measure on $I$ which is an eigenvector of $A^*$ for the
eigenvalue $\rho$, then $\nu=(\rho^{-n}\nu_0)$ is a state of $\underline
A$. By uniqueness of the state, $\nu_0=\mu_0$ and $\rho=\lambda$. This shows the
existence and the uniqueness of the Perron-Frobenius eigenvector and eigenvalue of
the transpose $A^*$. This also gives the same result for $A$ since $A$ is
primitive if and only if $A^*$ is so. Our proof of the Perron-Frobenius-Ruelle
theorem will follow the same pattern: we first establish that the state space of
the stationary system defined by the Ruelle operator $L_g:C(X)\rightarrow C(X)$
has a single element and then deduce the
existence and the uniqueness of the Perron-Frobenius eigenvector and eigenvalue of
the transpose $L_g^*$. 
\end{ex}

To conclude, let us quote some results on the affine representation of
dimension groups. They can be found in \cite{goo:poag} in a more general framework.
Given an arbitrary sequence $\underline B=(B_n)$ of Markovian
operators, we want to relate the inductive limit ${\cal E}={\cal E}(\underline
B)$ and the state space
$S=S(\underline B)$. As usual in the theory of compact convex sets, we define
$Aff(S)$ as the space of continuous real-valued affine functions on $S$. It is a
closed subspace of $C(S)$. We have the evaluation map $\theta:f\in {\cal E}\mapsto
\hat f\in Aff(S)$ defined by
$\hat f(\mu)=\mu(f)$ for $\mu\in S$. The elements of the kernel of $\theta$
are called {\it infinitesimals}.

\begin{prop}\label{affine representation} (see \cite{goo:poag}, Ch. 7)
With the  above notation,
\begin{enumerate}
\item The evaluation map $\theta$ has dense range in $Aff(S)$; moreover the image
of positive cone of ${\cal E}$ is dense in the positive cone of $Aff(S)$.
\item $Aff(S)$ is isomorphic as an ordered Banach space with unit to the
completion of
${\cal E}/Ker\theta$ with respect to the norm
$$\|f\|=\inf\{\|f_n\|: f_n\in C(X_n), j_n(f_n)=f\}.$$
\item  $(Aff(S),(\theta\circ j_n))$ is universal among the
$(F,(\varphi_n))$ where $F$ is an ordered Banach space with an order unit
$1$ and $\varphi_n:C(X_n)\rightarrow F$ are positive unital linear maps compatible
with the $B_n$'s which have the following property: if $f_k\in C(X_k)$
and
$g_l\in C(X_l)$ are such that
$$\|B_n\ldots B_{k+1}(f_k)-B_n\ldots B_{k+1}(g_l)\|$$ tends to
$0$ when $n$ tends to $\infty$, then $\varphi_k(f_k)=\varphi_l(g_l)$.
\end{enumerate}
\end{prop}

\begin{proof} The statements of $(i)$ are proved in \cite{goo:poag}, Theorem 7.9.

For $(ii)$, let $f\in {\cal E}$ and $n$ and $f_n\in C(X_n)$ such that $f=j_n(f_n)$. For
$\mu=(\mu_n)\in S$, we have $\hat f(\mu)=\mu_n(f_n)$, hence $\|hat f\|\le \|f\|$.
The equality $\|\hat f\|= \|f\|$ results from \cite{goo:poag}, Proposition 7.12
(e), and the equality of the norm $\|f\|_u$ defined there and our norm $\|f\|$.

For $(iii)$, we first observe that $(Aff(S),(\theta_n=\theta\circ j_n))$ has the
required properties. Let us just check the last property. Let
$f_k\in C(X_k)$ and
$g_l\in C(X_l)$ be as above. Let $\mu\in K$
and $n\ge k,l$. Then,
\begin{eqnarray}
|\theta_k(f_k)(\mu)-\theta_l(g_l)(\mu)|&=&|\mu(f_k\circ\pi_k-g_l\circ\pi_l)|\nonumber\\
&=& |(\pi_n)_*\mu(B_n\ldots B_{k+1}(f_k)-B_n\ldots B_{l+1}(g_l))|\nonumber\\
&\le&\|B_n\ldots B_{k+1}(f_k)-B_n\ldots B_{l+1}(g_l)\|\nonumber
\end{eqnarray}
Since this last quantity goes to $0$ when $n$ goes to $\infty$, we obtain 
$\theta_k(f_k)=\theta_l(g_l)$. Suppose now that $(F,(\varphi_n))$ satisfies also
the required properties. We have to define $\varphi:Aff(S)\rightarrow F$ such that
$\varphi_n=\varphi\circ\theta_n$. We show that $\varphi_0(f)$ depends only
on $\theta(f)$. Indeed, if $\theta(f)=\theta(g)$, then $f-g\in K_\perp$. According
to the lemma, $\|B_n\ldots B_1(f)-B_n\ldots B_1(g)\|$ goes to $0$. By hypothesis,
this implies that
$\varphi_0(f)=\varphi_0(g)$. Therefore, there is a well-defined map
$\varphi:Aff(S)\rightarrow F$ such that $\varphi\circ\theta(f)=\varphi_0(f)$. The
property
$\varphi_n=\varphi\circ\theta_n$ follows.
\end{proof}

Recall that a partially ordered abelian group $G$ is called {\it archimedean} if
$$x,y\in G\quad nx\le y\quad \forall n\in{\bf N}\Rightarrow x\le 0.$$ Note that the
ordered real vector spaces $C(X)$ and $Aff(S)$ are archimedean.

\begin{prop}\label{archimedean} (\cite{goo:poag}, Theorem 7.7) With the  above
notation, the following properties are equivalent
\begin{enumerate}
\item $0$ is the only infinitesimal element of ${\cal E}$.
\item ${\cal E}$ is archimedean.
\item The evaluation map $\theta:{\cal E}\rightarrow Aff(S)$ is injective. 
\end{enumerate}
\end{prop}


\begin{thebibliography}{10}

\bibitem{dr:amenable}  C.~ Anantharaman-Delaroche and J.~Renault: {\it Amenable
groupoids}, Monographie de l'Enseignement Math\'ematique No {\bf 36},
Gen\`eve, 2000.

\bibitem{bou:topo} N.~Bourbaki: {\it Topologie g\'en\'erale, chapitres 1\`a 4},
El\'ements de Math\'ematique, Diffusion C.C.L.S., Paris, 1971.

\bibitem{bou:Walters}  T.~Bousch: {\it La condition de Walters}, Janvier 2000

\bibitem{bd:g-measures}  G.~Brown and A.~H.~Dooley: {\it Odometer actions on
G-measures}, Ergodic Theory Dyn. System, {\bf 11} (1991), 279--307.

\bibitem{cap:gibbs}  D.~Capocaccia: {\it A definition of Gibbs state for a
compact set with ${\bf Z}^\nu$ action}, Commun. Math. Phys., {\bf 48}
(1976), 85--88.

\bibitem{cun:markovII}  J.~Cuntz: {\it A class of $C^*$-algebras and
topological Markov chains II: reducible chains and the Ext-functor for $C^*$-algebras},
Inventiones Math., {\bf  63} (1981), 25--40.

\bibitem{ckr:markovI}  J.~Cuntz and W.~Krieger: {\it A class of $C^*$-algebras and
topological Markov chains},
Inventiones Math., {\bf  56} (1980), 251--268.

\bibitem{dea:groupoids} V.~Deaconu: {\it Groupoids associated with
endomorphisms}, Trans. Amer. Math. Soc., {\bf 347} (1995), 1779--1786.

\bibitem{eff:dimensions}  E.~Effros: {\it Dimensions and C$^*$-algebras},
CBMS Regional Conf. Series in Math., {\bf 46}, Amer.
Math. Soc., 1981.

\bibitem{ell:ind} G.~Elliott: {\it On the classification of inductive limits of
sequences of semisimple finite-dimensional algebras}, J. Algebra, {\bf
38} (1976), 29--44.

\bibitem{efw:kms} M.~Enomoto, M.~Fujii and Y.~Watatani: {\it KMS states for
gauge action on $O_A$}, Math. Japon., {\bf 29} (1984), 607-619.

\bibitem{eva:on} D.~E.~Evans: {\it On $O_n$}, Publ. RIMS Kyoto Univ. {\bf 16} 
(1980), 915-927.

\bibitem{eva:markov} D.~E.~Evans: {\it The $C^*$-algebras of topological Markov
chains}, Tokyo Metropolitan University Lecture Notes (1983).

\bibitem{eva:quasiproduct} D.~E.~Evans: {\it Quasi-product states on
C$^*$-algebras}, Lecture Notes in Mathematics, No. 1132,
Springer-Verlag, Berlin-New York, 1985, 129-151. 

\bibitem{exe:RPF} R.~Exel: {\it KMS states for generalized gauge actions on
Cuntz-Krieger algebras}, preprint, (2001).

\bibitem{el:infinitematrices}  R.~Exel and M.~Laca: {\it Cuntz-Krieger algebras
for infinite matrices}, J. reine angew. Mth. (Crelle) {\bf 512} (1999), 119--172.

\bibitem{el:kms}  R.~Exel and M.~Laca: {\it Partial Dynamical Systems and the KMS
Condition}, preprint, (2000).

\bibitem{fam:sym}  T.~Fack and O.~Mar\'echal: {\it Sur la classification des
sym\'etrie des C$^*$-alg\`ebres UHF}, Canadian Journal of Math.

\bibitem{fel:proba}  W.~Feller: {\it An introduction to probability theory and
its applications}, vol. II, John Wiley, N.Y.

\bibitem{fan:Ruelle} A.~H.~Fan: {\it A proof of the Ruelle theorem}, Reviews Math.
Phys. {\bf 7}, no. 8 (1995), 1241--1247.

\bibitem{fj:Ruelle} A.~H.~Fan and Y.~P.~Jiang: {\it On Ruelle-Perron-Frobenius
Operators. I. Ruelle theorem}, Commun. Math. Phys. {\bf 223} (2001), 125--141.

\bibitem{geo:gibbs}  H.-O. Georgii: {\it Gibbs measures and phase transitions},
de Gruyter Studies in Mathematics, 1988.


\bibitem{gps:affable}  T.~Giordano, I.~Putnam and C.~Skau: {\it Affable equivalence relations and orbit structure of Cantor
dynamical systems}, preprint 2002.

\bibitem{goo:poag}  K.~R.~Goodearl: {\it Partially ordered abelian
groups with interpolation}, Mathematical Surveys and Monographs, {\bf 20},
Amer. Math. Soc., 1986.

\bibitem{kel:equilibrium}  G.~Keller: {\it Equilibrium States in Ergodic Theory},
London Mathematical Society Student Texts  {\bf 42}, Cambridge University Press
(1998).

\bibitem{ker:young}  S.~Kerov: {\it The boundary of Young lattice and random
Young tableaux}, DIMACS Ser. Discrete Math. Theoret. Comput. Sci.{\bf 24}, Amer.
Math. Soc., (1996), 133--158.

\bibitem{kp:pressure} D.~Kerr and C.~Pinzari: {\it Noncommutative pressure and the
variational principle in Cuntz-Krieger type C*-algebras}, J.~Funct. Anal., 
{\bf 188} (2002), 156-215.

\bibitem{kri:dimension}  W.~Krieger: {\it On dimension functions and topological
Markov chains},  Inventiones Math., {\bf  56} (1980), 239--250.

\bibitem{kum:localizations}  A.~Kumjian: {\it On localizations and simple
C*-algebras}, Pacific J. Math. {\bf 112} (1984), 141--192.

\bibitem{kr:kms}  A.~Kumjian and J.~Renault: {\it KMS states on C$^*$-algebras
associated to expansive maps}, in preparation.

\bibitem{lm:intro}  D.~Lind and B.~Marcus: {\it An introduction to symbolic
dynamics and coding}, Cambridge University Press, Cambridge, 1995.

\bibitem{mw:continuous trace}  P.~Muhly, D.~Williams: {\it Continuous trace
groupoid C*-algebras},  Math. Scand. {\bf 66} (1990), 231--241.

\bibitem{ren:approach}  J.~Renault: {\it A groupoid approach to
$C^*$-algebras}, Lecture Notes in Mathematics, Vol.~{\bf 793}
Springer-Verlag Berlin, Heidelberg, New York (1980).

\bibitem{ren:AF}  J.~Renault: {\it  AF equivalence relations and their
cocycles}, to appear in the Proceedings of the OAMP Conference (Constanta, 2001),
the Theta Foundation.

\bibitem{toe:AF}  A.~T\"or\"ok: {\it AF-algebras with unique trace}, Acta.
Sci. Math. (Szeged) {\bf 55}, No 1-2 (1991), 129--139.

\bibitem{vk:asymptotic}  A.~Vershik and S.~Kerov: {\it Asymptotic character
theory of the symmetric group}, Funct. Analysis and its Applications {\bf 15},
(1981), 25--36.

\bibitem{wal:78}  P.~Walters: {\it Invariant measures and equilibrium states for some
mappings which expand distances}, Trans. Amer. Math. Soc., {\bf 236}, (1978) ,
121--153.

\bibitem{wal:introduction}  P. Walters, {\it An introduction to ergodic theory},
Graduate Texts in Mathematics, {\bf 79}, Springer-Verlag, New York-Berlin, 1982. 

\bibitem{wal:Bowen}  P.~Walters: {\it Convergence of the Ruelle operator for a
function satisfying Bowen's condition}, Trans. Amer. Math. Soc., {\bf 353}, No 1
(2000) , 327--347.

\bibitem{wasa:thesis} A.~Wassermann: {\it Automorphic actions of compact groups
on operator algebras}, Ph.D thesis, U. of Pennsylvania (1981).

\bibitem{zin:thermo} M.~Zinsmeister: {\it Formalisme thermodynamique et
syst\`emes dynamiques holomorphes}, Panoramas et Synth\`eses, {\bf 4}, Soci\'et\'e
Math\'ematique de France, 1996.







\end{thebibliography}
\end{document}